\documentclass[a4paper, final, 12pt]{article}
\usepackage{amsmath}
\usepackage{amsfonts}
\usepackage{makeidx}
\usepackage{hyperref}
\usepackage{amsthm}
\usepackage{mathtools}
\usepackage{amssymb}
\usepackage{bm}
\usepackage{cite}
\usepackage{showkeys}
\usepackage[english]{babel}
\usepackage{dblaccnt}
\usepackage{accents}
\usepackage{graphicx}
\usepackage{psfrag}
\usepackage{subfig}
\usepackage{color}
\usepackage{float}
\usepackage{amscd}
\usepackage[mathscr]{euscript}
\usepackage{lipsum}
\usepackage{epstopdf}
\usepackage{algorithm}
\usepackage{algpseudocode}

\newtheorem{theorem}{Theorem}[section]
\newtheorem{lemma}{Lemma}[section]
\newtheorem{problem}{Problem}[section]
\newtheorem{assumption}{Assumption}[section]

\newtheorem{corollary}{Corollary}[section]
\newtheorem{remark}{Remark}[section]
\numberwithin{equation}{section}
\newcommand{\usc}{u_{\rm sc}}

\newcommand{\x}{{\bf x}}

\newcommand{\Pm}{P_{\rm meas}}
\newcommand{\Pp}{P_{\rm prop}}

\input{tcilatex}
\begin{document}

\title{A coefficient inverse problem with a single measurement of phaseless
scattering data}
\author{Michael V. Klibanov\thanks{
Department of Mathematics and Statistics, University of North Carolina at
Charlotte, Charlotte, NC 28223; (\texttt{mklibanv@uncc.edu}, \texttt{%
lnguye50@uncc.edu}) } \and Dinh-Liem Nguyen\thanks{%
Department of Mathematics, Kansas State University, Manhattan, KS 66506 (%
\texttt{dlnguyen@ksu.edu})} \and Loc H. Nguyen\footnotemark[1] }
\date{}
\maketitle

\begin{abstract}
This paper is concerned with a numerical method for a 3D coefficient inverse
problem with phaseless scattering data. These are multi-frequency data
generated by a single direction of the incident plane wave. Our numerical
procedure consists of two stages. The first stage aims to reconstruct the
(approximate) scattered field at the plane of measurements from its
intensity. We present an algorithm for the reconstruction process and prove
a uniqueness result of this reconstruction. After obtaining the approximate
scattered field, we exploit a newly developed globally convergent numerical
method to solve the coefficient inverse problem with the phased scattering
data. The latter is the second stage of our algorithm. Numerical examples
are presented to demonstrate the performance of our method. Finally, we
present a numerical study which aims to show that, under a certain
assumption, the solution of the scattering problem for the 3D scalar
Helmholtz equation can be used to approximate the component of the electric
field which was originally incident upon the medium.
\end{abstract}

%\begin{center}
%\bigskip \textbf{Liem, Loc,}
%\end{center}
%
%\textbf{There are a number of questions to you in the text in bold faced
%letters. Please address and then delete my questions.}
%
%1. \textbf{Proofs of Theorems 3.1, 3.2 and both Algorithms are encapsulated.
%I cannot work with encapsulated things. } \textcolor{red}{Done.}
%
%\textbf{2. Please do not give me in the future files with encapsulated
%parts. But this is for future papers only. \ As to this one, you can leave
%encapsulated parts here as they are. } \textcolor{red}{OK.}
%
%\textbf{3. On page 9 of the pdf file, after (3.21) \textquotedblleft from
%(3.19)-(5.6)" replace with \textquotedblleft from (3.19)-(3.21)"} \textcolor{red}{ Done.}
%
%\textbf{4. On page 10 \textquotedblleft By (3.20) and (3.22)" replace with
%\textquotedblleft By (3.20)-(3.22)"} \textcolor{red}{Done.}
%
%\textbf{5. Label for Figure 1: please do it in} \emph{italics}. \textcolor{red}{Done.}


\bigskip


\section{Introduction}

\label{sec: Intro}

The goal of this paper is to develop a new numerical method for a 3D
phaseless coefficient inverse scattering problem in the case when the data
to be inverted are generated by a single measurement event at multiple
frequencies. We assume that only the intensity, i.e. the square modulus, of
a complex valued wave field can be measured outside of scatterers and phase
cannot be measured. We use only a single direction of the incident plane
wave. In \ other words, we consider the phaseless coefficient inverse
scattering problem with single measurement data. Thus, this is a
non-overdetermined case, i.e. the number of free variables in the data
equals to the number of free variables in the unknown coefficient. We
propose a two-stage reconstruction procedure. In the first stage we
approximately reconstruct the scattered wave field at the plane of
measurements. Hence, this stage leads to a conventional phased coefficient
inverse scattering problem: when the whole complex valued wave field is
known at a part of the measurement plane. Next, to reconstruct the unknown
coefficient of Helmholtz equation, we apply the newly developed globally
convergent algorithm of \cite{Kliba2016}. According to \cite{BK,Kliba2016},
we call a numerical method for a coefficient inverse problem \emph{globally
convergent} if a theorem is proven, which guarantees that this method
delivers at least one point in a sufficiently small neighborhood of the
exact solution without any \emph{a priori} knowledge of this neighborhood.

Unlike the current paper, in \cite{KlibanovLocKejia:apnum2016} the case when
the intensity is given on an interval of frequencies for multiple point
sources was considered. While locations and shapes of unknown scatterers
were imaged accurately in \cite{KlibanovLocKejia:apnum2016}, the accuracy of
reconstructed abnormality/background contrasts was poor. This is because a
linearization of the travel time function was used in \cite%
{KlibanovLocKejia:apnum2016}. On the other hand, the globally convergent
numerical inversion method of \cite{Kliba2016} provides very accurate
locations and contrasts of abnormalities for a single measurement case. The
latter was consistently demonstrated on both computationally simulated \cite%
{Kliba2016} and experimental data \cite%
{AlekKlibanovLocLiemThanh:anm2017,LiemKlibanovLocAlekFiddyHui:jcp2017}
including the case when unknown targets were buried in a sandbox \cite%
{nguyen2017:iip2017}. The arguments in this paragraph are the reasons of our
choice of the two-stage procedure.

The study of the coefficient inverse problem with phaseless scattering data
is motivated by applications in, e.g. imaging of nano-scale structures and
biological cells. Typical nano structures of interest have sizes of hundreds
of nanometers ($nm$). Recall that $1$ micron ($\mu m$)$=10^{3}nm.$ Typical
sizes of biological cells are in the range of $\left( 5,100\right) \mu m$ 
\cite{PM,Bio}. To image these, one should use optical sources with the same
range of wavelengths. However, the corresponding frequency is very large.
For example, the wavelength $\lambda =1\mu m$ corresponds to the frequency $%
\omega =299,792$ GHz. Hence, only the intensity of the scattered field can
be measured while the phase is lost \cite{Darahanau:PL2005, Dierolf:EN2008,
Khachaturov:cmmp2009, Petersena:U2008, Ruhlandt:PR2014}.

Solution of the coefficient inverse scattering problem without the phase
information is a long standing problem. For the first time, this problem was
probably posed in \cite[Chapter 10]{Chadan:sv1977}. The first uniqueness
result for this problem was established in \cite{Klibanov:jmp1992} in the 1D
case, also see \cite{AktosunSacks:ip1998} for a follow up result. In 3D, the
first uniqueness result was obtained in \cite{Klibanov:SiamJMath2014}.
Later, uniqueness theorems in 3D were established in \cite%
{Klibanov:AppliedMathLetters2014,Klibanov:aa2014,Klibanov:ipi2017,KlibRom,Zhang}%
. The analytic reconstruction procedures in 3D were proposed in \cite%
{KR:JIpp2015,
KlibanovRomanov:SIAMam2016,KlibanovRomanov:ejmca2015,KlibanovRomanov:ip2016}%
. However, these procedures require a large range of frequencies which might
be unrealistic in practice. Hence, the method of \cite%
{KlibanovRomanov:SIAMam2016} was modified and made suitable for computations
for a realistic range of parameters in \cite{KlibanovLocKejia:apnum2016}.

In publications \cite{Nov3,Nov4,Nov5} phaseless coefficient inverse
scattering problems were considered. Their statements are different from the
ones in papers cited above. Uniqueness theorems were proved and
reconstruction procedures were proposed in \cite{Nov3,Nov4,Nov5}. In \cite%
{BardsleyVasquez:ip2016, BardsleyVasquez:sjis2016} a phaseless coefficient
inverse scattering problem for Helmholtz equation was solved numerically
using Kirchhoff migration and Born approximation. While coefficients of PDEs
are subjects of interests in the above cited works, there is also a
significant interest in the reconstruction of surfaces of scatterers from
phaseless data. In this regard we refer to, e.g. publications \cite%
{AmmariChowZou:sjap2016,Bao,Iv1,Iv2,LLW,Zhang1} and references cited therein.

Let $k>0$ be the wave number. In fact, on the first stage of our procedure
we reconstruct the first term of the asymptotic expansion at $k\rightarrow
\infty $ of the solution of Helmholtz equation at the measurement plane. We
prove a theorem which claims that this reconstruction is unique. In this
paper, for simplicity, we use the Helmholtz equation to model the light
propagation. On the other hand, it is well-known that the wave field is
governed by the Maxwell's system. To validate our analysis, we numerically
compare in Appendix the solution of the Helmholtz equation with the one of
the Maxwell's system. Another validation can be attributed to quite accurate
results obtained by this group for phased microwave experimental data \cite%
{AlekKlibanovLocLiemThanh:anm2017,LiemKlibanovLocAlekFiddyHui:jcp2017,nguyen2017:iip2017}
in which the Helmholtz equation was used and the globally convergent
numerical method of \cite{Kliba2016} was applied.

In the next section, we formulate the phaseless coefficient inverse
scattering problem. In Section \ref{sec: uniqueness}, we prove a uniqueness
result of the phase retrieval. In Section \ref{sec: phase retrieval}, we
describe our numerical approach to reconstruct the lost phase. In Section %
\ref{sect: phased ICP}, we briefly summarize our globally convergent method
of \cite{Kliba2016} for the reader's convenience. In Section \ref{sect: Num}%
, we present our numerical results. In Section \ref{sec: 8} we provide
summary of our results. Finally, in the Appendix, which is Section \ref{sec:
Maxwell}, we compare numerically solutions to Maxwell's system and Helmholtz
equation.

\section{Problem Statement}

\label{sec: problem statement}

%
%In this paper, we are interested in the problem of identifying the nano-structures by consider them as scatterers and then quantitatively reconstructing their spatially distributed dielectric constants from the measurement of the intensities of the scattered lights.  
%The physical experiment for the purpose is described as follows.

\begin{figure}[h]
\begin{center}
\includegraphics[width = 0.5\textwidth]{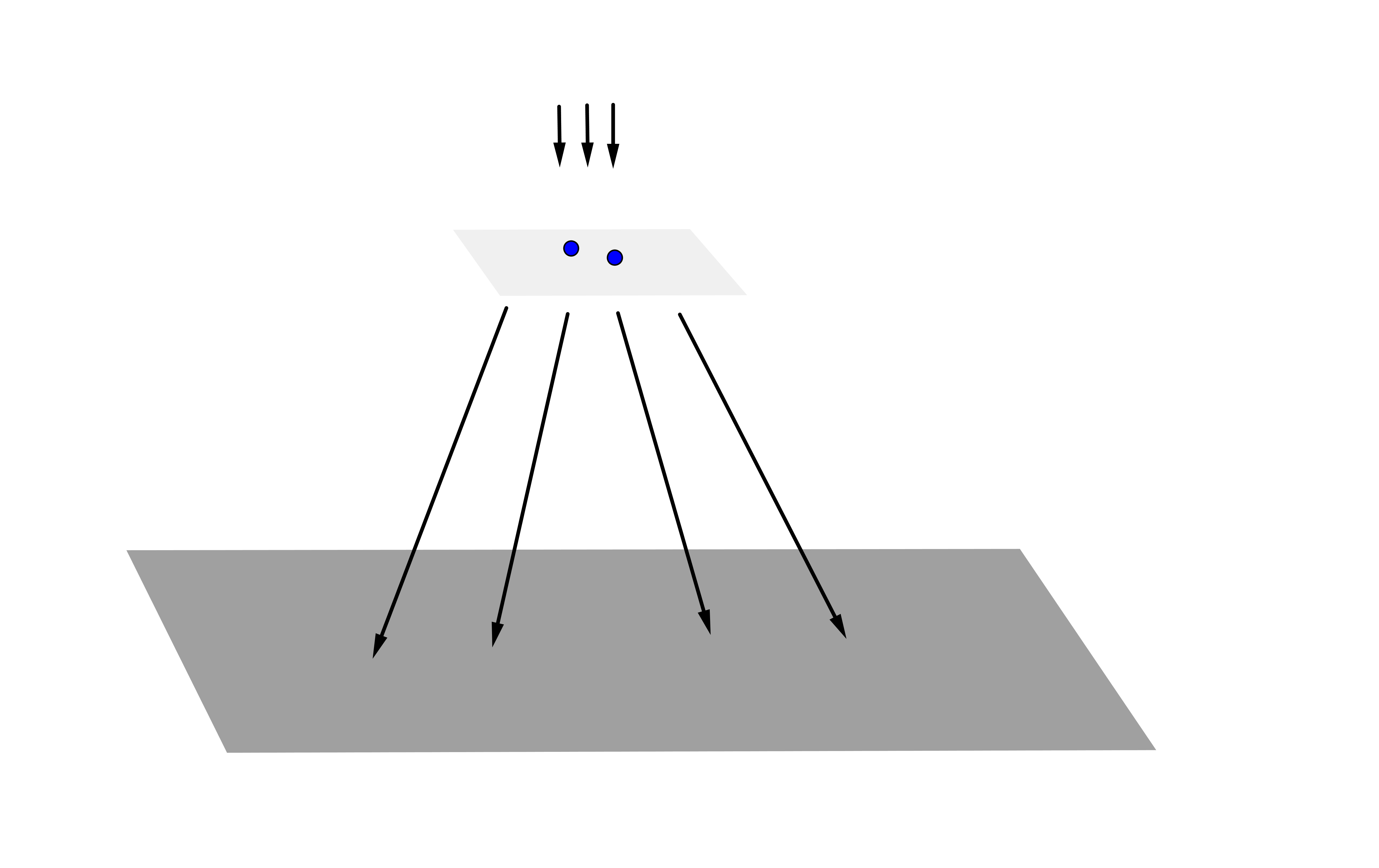} \put(-70,25){$\Pm$}
\put(-100,70){Scattered waves}
\end{center}
\caption{\it The experimental setup. The laser beam hits scatterer (dots) and
causes scattered wave. The intensity is measured on $P_{\mathrm{meas}}$. }
\label{fig 1}
\end{figure}

Let a laser beam illuminate the unknown nano-structure/biological cell,
which plays the role of a scatterer. The diameter of a laser beam is a few
millimeters (mm) and $1mm=10^{3}\mu m.$ Given that sizes of our scatterers
do not exceed 100 $\mu m=0.1mm$ (section 1), these scatterers
\textquotedblleft percept" that laser beam as a perfect plane wave. The
laser beam scatters after hitting the scatterers. It is well known that
modern light detectors placed inside of the laser beam are burned. Hence,
one should place detectors outside of that beam. However, outside of the
laser beam the total wave field approximately equals the scattered wave
field. Hence, we assume below that we measure the intensity of the scattered
wave at a square $P_{\mathrm{meas}}$ of a fixed plane, see Figure \ref{fig 1}
for an illustration. The problem we consider in this paper is to reconstruct
the spatially distributed dielectric constant of the scatterer from the
measurements of the intensities of the scattered waves on $P_{\mathrm{meas}}$
at an interval of frequencies. We point out, however, that a precise
mathematical modeling of the laser beam is outside of the scope of this
publication. So, the above considerations were given only to explain why do
we consider the intensity of the scattered rather than of the total wave
field. 

\subsection{The phaseless coefficient inverse scattering problem}

\label{subsec: inverse problem statement}

Denote $\mathbf{x=}(x_{1},x_{2},x_{3})\in \mathbb{R}^{3}.$ Let $\Omega $ be
a bounded domain in $\mathbb{R}^{3}$ with a smooth boundary $\partial \Omega 
$ and such that $\Omega \subset \left\{ x_{3}>0\right\} .$ Let $c(\mathbf{x}%
),$ be a function satisfying the following conditions: 
\begin{equation}
c(\mathbf{x})\in C^{15}(\mathbb{R}^{3}),\quad c(\mathbf{x})\geq 1%
\mbox{
for all }\mathbf{x}\in \mathbb{R}^{3}\quad \mbox{ and }c(\mathbf{x})=1%
\mbox{
for all }\mathbf{x}\in \mathbb{R}^{3}\setminus \Omega .  \label{condition c}
\end{equation}%
The function $c(\mathbf{x})$ models the spatially distributed dielectric
constant of the medium consisting of the background and the scatterers. The
condition $c(\mathbf{x})=1$ in $\mathbb{R}^{3}\setminus \Omega $ means that
we suitably scale the dielectric constant, so that it equals to 1 in the
background (vacuum). The condition $c(\mathbf{x})\geq 1$ means that the
dielectric constant in the medium is not less than in the vacuum. The
Riemannian metric corresponding to the function $c(\mathbf{x})$ is given by 
\begin{equation*}
d\tau (\mathbf{x})=\sqrt{c(\mathbf{x})}|d\mathbf{x}|,\quad |d\mathbf{x}|=%
\sqrt{(dx_{1})^{2}+(dx_{2})^{2}+(dx_{3})^{2}}.
\end{equation*}%
Fix the number $a>0.$ Consider the plane $P_{a}=%
\{(x_{1},x_{2},-a):x_{1},x_{2}\in \mathbb{R}\}.$We impose the following
condition on the function $c(\mathbf{x})$: 
%Let $a$ be a arbitrary positive number. Consider the plane 
%\[
%	P_{-a} = \{\x \in \R^3: x_3 = -a\}.
%\] 
%Since $c(x) = 1$ outside $\Omega$, the geodesic lines are straight in $\R^3 \setminus \overline \Omega$. In particular, since our incident planewave propagates along the positive direction of the $x_3-$axis, geodesic lines which are parallel to the $x_3-$axis in the half space $\{x_3 < 0\}$

\begin{assumption}[Assumption of Regularity of Geodesic Lines.]
\label{assumption regularity} For any point $\mathbf{x}\in \mathbb{R}^{3}$
there exists a unique geodesic line $\Gamma (\mathbf{x},a)$, with respect to
the metric $d\tau $, connecting $\mathbf{x}$ with the plane $P_{a}$ and
perpendicular to $P_{a}$.
\end{assumption}

The following sufficient condition of the regularity of geodesic lines was
derived in \cite{Rom3} 
\begin{equation*}
\dsum\limits_{i,j=1}^{3}\frac{\partial ^{2}c\left( \mathbf{x}\right) }{%
\partial x_{i}\partial x_{j}}\xi _{i}\xi _{j}\geq 0,\forall \mathbf{x}\in 
\overline{\Omega },\forall \mathbf{\xi }\in \mathbb{R}^{3}.
\end{equation*}

\begin{remark}
{The smoothness condition (\ref{condition c}) imposed on the function }$c(%
\mathbf{x})$ {as well as Assumption \ref{assumption regularity} are
necessary for the theoretical purposes only: to derive the asymptotic
behavior (\ref{4}). However, we do not verify neither condition (\ref%
{condition c}) nor Assumption \ref{assumption regularity} in our numerical
studies. Indeed, that asymptotic behavior is derived in Theorem 3.1 on the
basis of the construction of the solution of the Cauchy problem for a
hyperbolic equation. However, this construction requires that }$c(\mathbf{x}%
)\in C^{15}(\mathbb{R}^{3})$ {\cite{KlibanovRomanov:ip2016,Rom2}. Besides,
the minimal smoothness of unknown coefficients is usually not of a great
concern in studies of coefficient inverse problems, see, e.g. \cite%
{Nov1,Nov2} and theorem 4.1 in \cite{Rom1}.}
\end{remark}

\bigskip The function $\tau (\mathbf{x})$\ is the travel time from the plane 
$P_{a}$ to the point $\mathbf{x}$ and \cite{KlibanovRomanov:ip2016}\emph{\ }%
\begin{equation}
\tau (\mathbf{x})=\dint\limits_{\Gamma (\mathbf{x},a)}\sqrt{c\left( \mathbf{%
\xi }\right) }d\sigma .  \label{200}
\end{equation}

\begin{lemma}
{For all }$\mathbf{x}\in \left\{ x_{3}>-a\right\} $ we have $\tau (\mathbf{x}%
)\geq x_{3}.${\ Consider the set }$\Omega _{1}, $ 
\begin{equation}
\Omega _{1}=\left\{ \mathbf{x}:c\left( \mathbf{x}\right) >1\right\} \subset
\Omega .  \label{201}
\end{equation}%
{Assume that the set (\ref{201}) is convex and its boundary} $\partial
\Omega _{1}\in C^{1}.$ {Then }$\tau (\mathbf{x})=x_{3}$ {\ if and only if} $%
\Gamma \left( \mathbf{x},a\right) =L\left( \mathbf{x},a\right) ,$ {where} $%
L\left( \mathbf{x},a\right) $ {is the straight line connecting the point }$%
\mathbf{x}$ {\ with the plane }$P_{a}$ {and orthogonal to }$P_{a}.$ If $\tau
(\mathbf{x})=x_{3},$ {then} {for any point }$\mathbf{x}^{\prime }\mathbf{\in 
}\Gamma \left( \mathbf{x},a\right) =L\left( \mathbf{x},a\right) ,$ {except
of probably one point, there exists such a sufficiently small neighborhood} $%
O\left( \mathbf{x}^{\prime }\right) $ {of} $\mathbf{x}^{\prime }$ {that} $%
\tau (\mathbf{x}^{\prime \prime })=x_{3}^{\prime \prime },\forall \mathbf{x}%
^{\prime \prime }=\left( x_{1}^{\prime \prime },x_{2}^{\prime \prime
},x_{3}^{\prime \prime }\right) \in O\left( \mathbf{x}^{\prime }\right) .$
\end{lemma}

%\begin{proof}
\noindent {\it Proof.}
 The inequality $\tau (\mathbf{x})\geq x_{3}$ follows from (\ref%
{200}) as well as from the fact that by (\ref{condition c}) $c\left( \mathbf{%
x}\right) \geq 1$. Suppose now that $\tau (\mathbf{x})=x_{3}.$ If $\Gamma
\left( \mathbf{x},a\right) \cap \Omega _{1}\neq \varnothing ,$ then (\ref%
{200}) implies that $\tau (\mathbf{x})>x_{3}.$ Hence, $\Gamma \left( \mathbf{%
x},a\right) \cap \Omega _{1}=\varnothing .$ This means that $\Gamma \left( 
\mathbf{x},a\right) =L\left( \mathbf{x},a\right) .$ Now, either $L\left( 
\mathbf{x},a\right) \cap \partial \Omega _{1}=\varnothing $ or, due to both
the convexity of the domain $\Omega _{1}$ and the smoothness of its boundary 
$\partial \Omega _{1},$ the straight line $L\left( \mathbf{x},a\right) $ is
the tangent line to $\partial \Omega _{1}$ at a certain unique point.
Obviously in both these two cases the assertion of this lemma about $O\left( 
\mathbf{x}^{\prime }\right) $ is true. $\square$
%\end{proof}
Let $[\underline{k},\overline{k}]$ be an interval of wave numbers 
$k=2\pi /\lambda \in \lbrack \underline{k},\overline{k}]$ where $\lambda $
is the dimensionless wavelength. Consider the incident plane wave $u_{%
\mathrm{inc}}(\mathbf{x},k)$ propagating along the $x_{3}$ axis, 
\begin{equation}
u_{\mathrm{inc}}(\mathbf{x},k)=\exp (\mathrm{i}kx_{3}).  \label{eqn source}
\end{equation}%
The propagation of the total wave field $u(\mathbf{x},k)$ is governed by the
Helmholtz equation and the outgoing Sommerfeld radiation condition, 
\begin{equation}
u(\mathbf{x},k)=u_{\mathrm{inc}}(\mathbf{x},k)+u_{\mathrm{sc}}(\mathbf{x}%
,k),\quad \mathbf{x}\in \mathbb{R}^{3},k\in \lbrack \underline{k},\overline{k%
}],  \label{eqn total field}
\end{equation}%
\begin{equation}
\left\{ 
\begin{array}{rcll}
\Delta u(\mathbf{x},k)+k^{2}c(\mathbf{x})u(\mathbf{x},k) & = & 0, & \mathbf{x%
}\in \mathbb{R}^{3}, \\ 
\partial _{r}u_{\mathrm{sc}}(\mathbf{x},k)-\mathrm{i}ku_{\mathrm{sc}}(%
\mathbf{x},k) & = & o(r^{-1}), & \mbox{as }r=|\mathbf{x}|\rightarrow \infty .%
\end{array}%
\right.  \label{eqn Helmholtz}
\end{equation}%
Let the number $R>0.$ Denote 
\begin{equation}
P_{\mathrm{meas}}=\{\mathbf{x}=(x_{1},x_{2},x_{3}):-b<x_{1},x_{2}<b,x_{3}=R\}
\label{measurement plane}
\end{equation}%
the square on the plane $P=\left\{ x_{3}=R\right\} $ where measurements of
the intensity $\left\vert u(\mathbf{x},k)\right\vert ^{2}$ are conducted.
Here $R>0$ is the distance from the origin to the measurement plane and $b>0$
is the size of that rectangle. Assume that the plane $P$ does not intersect
with $\Omega $, $\overline{\Omega }\cap P=\varnothing .$ The phaseless
coefficient inverse scattering problem is formulated as:

\textbf{Phaseless Coefficient Inverse Scattering Problem (PCISP)}. \emph{%
Given the data} 
\begin{equation}
f(\mathbf{x},k)=|u_{\mathrm{sc}}(\mathbf{x},k)|^{2},\quad \mathbf{x}\in P_{%
\mathrm{meas}},k\in \lbrack \underline{k},\overline{k}],  \label{eqn f}
\end{equation}%
\emph{determine the dielectric constant} $c(\mathbf{x})$ \emph{for} $\mathbf{%
x}\in \Omega .$

\begin{remark}
\begin{enumerate}
\item {We model the wave propagation by the single Helmholtz equation with
the outgoing radiation condition instead of the full Maxwell's system, see
section 9 and a discussion in Introduction.}

\item {A natural question about the uniqueness of the PCISP arises. We prove
in Section \ref{sec: uniqueness} that one can uniquely reconstruct the first
term of the asymptotic expansion at }$k\rightarrow \infty ${\ of the
function }$u(\mathbf{x},k),\mathbf{x}\in P_{\mathrm{meas}}${\ from the data (%
\ref{eqn f}). Let the square }$P_{\mathrm{meas}}\subset P,$ where $P$ is the
corresponding plane. {The next question, however, is about the uniqueness of
the reconstruction of the coefficient }$c(\mathbf{x})${\ for }$\mathbf{x}\in
\Omega ${\ even in the case when the whole function }$u(\mathbf{x},k)${\
(rather than that first term only) is known for all }$\mathbf{x}\in P${\ and
for all }$k>0.${\ Addressing this question is a well known long standing
open problem. Indeed, all uniqueness theorems for }$n-D${, }$n\geq 2${\
coefficient inverse problems with single measurement data are currently
proven only by the method, which was originally proposed in \cite{BukhKlib}
in 1981, also see, e.g. the section 1.10 in the book \cite{BK}, the book 
\cite{KT}, the survey \cite{Ksurvey} and references cited in \cite{Ksurvey}.
Carleman estimates are the key ingredient of this method. However, in our
specific case, this method works only if the right hand side of Helmholtz
equation (\ref{eqn Helmholtz}) is non vanishing in }$\overline{\Omega }.${\
Hence, we just assume uniqueness of that second problem: for the purpose of
computations. }
\end{enumerate}
\end{remark}

\subsection{The Lippman-Schwinger equation}

Assume that the function $c(\mathbf{x})$ satisfying \eqref{condition c} is
known for all $\mathbf{x}\in \mathbb{R}^{3}$. It follows from the Helmholtz
equation in \eqref{eqn Helmholtz} that%
\begin{equation*}
\Delta u_{\mathrm{sc}}(\mathbf{x},k)+k^{2}u_{\mathrm{sc}}(\mathbf{x}%
,k)+k^{2}\beta (\mathbf{x})u(\mathbf{x},k)=0,\forall \mathbf{x}\in \mathbb{R}%
^{3},
\end{equation*}%
\begin{equation}
\beta (\mathbf{x})=c(\mathbf{x})-1.  \label{1000}
\end{equation}%
Here, we have used \eqref{eqn total
field} and the fact that $\Delta u_{0}(\mathbf{x},k)+k^{2}u_{0}(\mathbf{x}%
,k)=0$. This and the outgoing Sommerfeld radiation condition in 
\eqref{eqn
Helmholtz} imply that 
\begin{equation}
u_{\mathrm{sc}}(\mathbf{x},k)=k^{2}\dint\limits_{\mathbb{R}^{3}}\frac{\exp (%
\mathrm{i}k|\mathbf{x}-\mathbf{\xi }|)}{4\pi |\mathbf{x}-\mathbf{\xi }|}%
\beta (\mathbf{x})u(\mathbf{\xi },k)d\mathbf{\xi .}  \label{eqn LS usc}
\end{equation}%
Using \eqref{eqn total field} and \eqref{eqn LS usc}, we obtain the
Lippmann-Schwinger equation 
\begin{equation}
u(\mathbf{x},k)=u_{\mathrm{inc}}(\mathbf{x},k)+k^{2}\dint\limits_{\Omega }%
\frac{\exp (\mathrm{i}k|\mathbf{x}-\mathbf{\xi }|)}{4\pi |\mathbf{x}-\mathbf{%
\xi }|}\beta (\mathbf{\xi })u(\mathbf{\xi },k)d\mathbf{\xi }.  \label{eqn LS}
\end{equation}%
For an integer $k\geq 0$ and for $\alpha \in \left( 0,1\right) $ let $%
C^{k+\alpha }(\mathbb{R}^{3})$ be H\"{o}lder spaces. The following result
holds \cite[Chapter 8]{Colto2013}:

\begin{theorem}
Assume that $c(\mathbf{x})\in C^{\alpha }(\mathbb{R}^{3})$ satisfies the
rest of conditions \eqref{condition c}. Then the Lippmann-Schwinger equation %
\eqref{eqn LS} has unique solution $u(\mathbf{x},k)\in C^{2+\alpha }(\mathbb{%
R}^{3})$ for all $k>0$. Moreover, this function $u(\mathbf{x},k)$ is the
unique solution of the problem \eqref{eqn Helmholtz}.
\end{theorem}

When using the globally convergent method \cite{Kliba2016} below, we solve
equation \eqref{eqn LS} on each iteration. To solve integral equation %
\eqref{eqn LS} numerically, we use the method developed in \cite{Lechl2014,
Nguye2014a}.

\section{Uniqueness result}

\label{sec: uniqueness}

In this section, we prove that the first term of the asymptotic expansion of
the function $u\left( \mathbf{x},k\right) $ at $k\rightarrow \infty $ for
points $\mathbf{x}\in P_{\text{meas}}$ can be determined uniquely from the
data $f(\mathbf{x},k).$ For any number $\theta >0$ we define the half-plane $%
\mathbb{C}_{\theta }$ of the complex plane $\mathbb{C}$ as%
\begin{equation*}
\mathbb{C}_{\theta }=\left\{ z\in \mathbb{C}:\Im z>-\theta \right\} .
\end{equation*}

\begin{theorem}
Assume that the function $c\left( \mathbf{x}\right) $ satisfies conditions %
\eqref{condition c}. Suppose that the Assumption of the Regularity of
Geodesic Lines holds true. Let $G\subset \{\mathbf{x}%
=(x_{1},x_{2},x_{3}):x_{3}>-a\}$ be an arbitrary bounded domain such that $%
P_{\text{meas}}\subset G$. Then there exists a number $\theta =\theta (G)>0$
such that for all points $\mathbf{x}\in G$, the solution $u(\mathbf{x},k)$
of problem \eqref{eqn total field}--\eqref{eqn Helmholtz} is analytic with
respect to $k\in \mathbb{R}$ and can be analytically continued in the
half-plane $\mathbb{C}_{\theta }.$ Furthermore, the function $f(\mathbf{x}%
,k)=|u_{sc}(\mathbf{x},k)|^{2},\mathbf{x}\in P_{\text{meas}}$ is analytic
with respect to $k\in \mathbb{R}.$ In addition, the following asymptotic
behavior holds 
\begin{equation}
u_{\mathrm{sc}}(\mathbf{x},k)=A(\mathbf{x})e^{\mathrm{i}k\tau (\mathbf{x}%
)}-e^{\mathrm{i}kx_{3}}+\mu (\mathbf{x},k),\text{ \ }\mathbf{x}\in
G,k\rightarrow \infty ,  \label{4}
\end{equation}%
where the function $A(\mathbf{x})>0$, the function $\tau (\mathbf{x})$ is
defined in Subsection \ref{subsec: inverse problem statement}, and the
function $\mu (\mathbf{x},k)$ is such that for $j=0,1,2$ 
\begin{equation}
\partial _{k}^{j}\mu (\mathbf{x},k)=O(k^{-1}),\quad \text{\ }\mathbf{x}\in
G,k\rightarrow \infty .  \label{5}
\end{equation}%
\label{thm 1}
\end{theorem}

In the proof of this theorem, we modify the material of section 4 of \cite%
{KlibanovRomanov:ip2016}. Although properties (\ref{4}), (\ref{5}) and the
analyticity of the function $u(\mathbf{x},k)$ follow from results of \cite%
{KlibanovRomanov:ip2016}, they are not explicitly formulated there.

%\begin{proof}[Proof of Theorem \ref{thm 1}]

\noindent{\it Proof of Theorem \ref{thm 1}.}
Consider the following auxiliary hyperbolic equation 
\begin{equation}
c(\mathbf{x})v_{tt}=\Delta v,\quad \mathbf{x}\in \mathbb{R}^{3},t\in \mathbb{%
R}.  \label{6}
\end{equation}%
And consider the solution of equation (\ref{4}) in the form 
\begin{equation}
v(\mathbf{x},t)=\delta \left( t-x_{3}\right) +\widetilde{v}(\mathbf{x},t),
\label{7}
\end{equation}%
where $\widetilde{v}(\mathbf{x},t)$ is such that 
\begin{equation}
\widetilde{v}\left( \mathbf{x},t\right) =0\quad \text{ for }x_{3}<-a.
\label{8}
\end{equation}%
Let $T>0$ be an arbitrary number. Denote%
\begin{equation*}
D\left( T\right) =\left\{ \left( \mathbf{x},t\right) :\max \left( -a,\tau
\left( \mathbf{x}\right) \right) <t<T\right\} ,
\end{equation*}%
\begin{equation*}
H\left( t\right) =\left\{ 
\begin{array}{ll}
1 & t>0, \\ 
0 & t<0.%
\end{array}%
\right.
\end{equation*}%
It was proven in \cite[Theorem 1]{KlibanovRomanov:ip2016} that the problem (%
\ref{6})-(\ref{8}) has unique solution of the form%
\begin{equation}
v\left( \mathbf{x},t\right) =A\left( \mathbf{x}\right) \delta \left( t-\tau
\left( \mathbf{x}\right) \right) +H\left( t-\tau \left( \mathbf{x}\right)
\right) \widehat{v}\left( \mathbf{x},t\right) ,  \label{9}
\end{equation}%
where the function $A\left( \mathbf{x}\right) >0$ and%
\begin{equation}
\widehat{v}\left( \mathbf{x},t\right) \in C^{2}\left( \overline{D\left(
T\right) }\right) .  \label{90}
\end{equation}

Furthermore, Theorem 4 of Chapter 10 of \cite{Vainberg:gbsp1989} as well as
Remark 3 after that theorem guarantee that there exists a number $\theta
=\theta \left( c,G\right) >0$ and a number $C_{1}=C_{1}\left( c,G\right)
>0,C_{2}=C_{2}\left( c,G\right) >0,$ all three depending only on listed
parameters, such that%
\begin{equation}
\left\vert D_{\mathbf{x}}^{\alpha }D_{t}^{k}\widehat{v}\left( \mathbf{x}%
,t\right) \right\vert \leq C_{2}e^{-\theta t},\quad \forall \mathbf{x}\in
G,\forall t>C_{1},\left\vert \alpha \right\vert +k\leq 2.  \label{10}
\end{equation}%
Here $\alpha =\left( \alpha _{1},\alpha _{2},\alpha _{3}\right) $ is the
multiindex with non-negative integer coordinates and $\left\vert \alpha
\right\vert =\alpha _{1}+\alpha _{2}+\alpha _{3}.$ By (\ref{9}) and (\ref{10}%
) we can consider Fourier transform of the function $v\left( \mathbf{x}%
,t\right) ,$%
\begin{equation}
V\left( \mathbf{x},k\right) =\dint\limits_{-\infty }^{\infty }v\left( 
\mathbf{x},t\right) e^{\mathrm{i}kt}dt=A\left( \mathbf{x}\right) e^{\mathrm{i%
}k\tau \left( \mathbf{x}\right) }+\dint\limits_{\tau \left( \mathbf{x}%
\right) }^{\infty }\widehat{v}\left( \mathbf{x},t\right) e^{\mathrm{i}kt}dt.
\label{11}
\end{equation}%
Next, Theorem 3.3 of \cite{Vainberg:rms1966} and Theorem 6 of Chapter 9 of 
\cite{Vainberg:gbsp1989} imply that 
\begin{equation}
V\left( \mathbf{x},k\right) =u\left( \mathbf{x},k\right) ,\quad \forall 
\mathbf{x}\in \mathbb{R}^{3},\forall k>0,  \label{12}
\end{equation}%
where $u\left( \mathbf{x},k\right) $ is the solution of our original forward
problem \eqref{eqn total field}--\eqref{eqn Helmholtz}. Next, using (\ref{9}%
)-(\ref{10}), (\ref{12}) and the integration by parts in (\ref{11}), we
obtain 
\begin{equation*}
u\left( \mathbf{x},k\right) =A\left( \mathbf{x}\right) e^{\mathrm{i}k\tau
\left( \mathbf{x}\right) }+\frac{\mathrm{i}}{k}\widehat{v}\left( \mathbf{x}%
,\tau \left( \mathbf{x}\right) \right) e^{\mathrm{i}k\tau \left( \mathbf{x}%
\right) }+\frac{\mathrm{i}}{k}\dint\limits_{\tau \left( \mathbf{x}\right)
}^{\infty }\partial _{t}\widehat{v}\left( \mathbf{x},t\right) e^{\mathrm{i}%
kt}dt,\quad k>0,
\end{equation*}%
which proves the asymptotic expansion (\ref{4}), (\ref{5}).

Next, it follows from (\ref{9})--(\ref{12}) that%
\begin{equation*}
\partial _{k}V\left( \mathbf{x},k\right) =\mathrm{i}A\left( \mathbf{x}%
\right) \tau \left( \mathbf{x}\right) e^{\text{i}k\tau \left( \mathbf{x}%
\right) }+\mathrm{i}\dint\limits_{\tau \left( \mathbf{x}\right) }^{\infty }%
\widehat{v}\left( \mathbf{x},t\right) te^{\mathrm{i}kt}dt,\quad \forall k\in 
\mathbb{C}_{\theta }.
\end{equation*}%
Hence, using (\ref{12}), we conclude that the function $u\left( \mathbf{x}%
,k\right) $ has analytic continuation with respect to $k$ from the real line 
$\mathbb{R}$ in the half-plane $\mathbb{C}_{\theta }.$ Finally, the
analyticity of the function $f\left( \mathbf{x},k\right) =\left\vert u_{%
\mathrm{sc}}\left( \mathbf{x},k\right) \right\vert ^{2}$ with respect to $%
k\in \mathbb{R}$ follows from \cite[Lemma 3.5]{Klibanov:ipi2017}. 
$\square$
%\end{proof}

Theorem 3.2 is similar with theorem 1 of \cite[Lemma 3.5]{KlibRom}. While
theorem 1 of \cite[Lemma 3.5]{KlibRom} works for the case when the wave
field is generated by a point source, Theorem 3.2 is valid for the case of
the incident plane wave. The proof here is different from the one in \cite[%
Lemma 3.5]{KlibRom}. Two major differences are that neither an analog of
Lemma 2.1 nor the derivative $\varphi ^{\prime }\left( k\right) $ in (\ref%
{17}) were not considered in \cite[Lemma 3.5]{KlibRom}.

\begin{theorem}
Assume that the function $c\left( \mathbf{x}\right) $ satisfies condition %
\eqref{condition c}. Suppose that the Assumption of the Regularity of
Geodesic Lines holds true. In addition, assume that conditions of Lemma 2.1
hold. Consider an arbitrary point $\mathbf{x}_{0}\in P_{\text{meas}}.$ Also,
consider the function $\varphi \left( k\right) ,$ 
\begin{equation}
\varphi \left( k\right) =\left\vert u_{sc}\left( \mathbf{x}_{0},k\right)
\right\vert ^{2},k\in \left( \underline{k},\overline{k}\right) ,  \label{13}
\end{equation}%
where $\left( \underline{k},\overline{k}\right) \subset \mathbb{R}$ is a
certain interval. Then the numbers $\tau \left( \mathbf{x}_{0}\right) $ and $%
A\left( \mathbf{x}_{0}\right) $ in the asymptotic expansion \eqref{4}, %
\eqref{5} are uniquely determined from the knowledge of the function $%
\varphi \left( k\right) $ in \eqref{13}. Furthermore, if $\tau \left( 
\mathbf{x}_{0}\right) =x_{3,0}$ for $\mathbf{x}_{0}=\left(
x_{1,0},x_{2,0},x_{3,0}\right) ,$ then $A\left( \mathbf{x}_{0}\right) =1.$ %
\label{thm 2}
\end{theorem}

Corollary \ref{col 1} follows immediately from Theorem \ref{thm 2}.

\begin{corollary}
\label{col 1} Consider the PCISP. Then functions $A\left( \mathbf{x}\right) $
and $\tau \left( \mathbf{x}\right) $ in the asymptotic expansion \eqref{4}, %
\eqref{5} are uniquely determined for $\mathbf{x}\in P_{\mathrm{meas}}$ from
the knowledge of the function $f\left( \mathbf{x},k\right) $ in \eqref{eqn f}%
.
\end{corollary}

%\begin{proof}{Proof of Theorem \ref{thm 2}}
\emph{Proof of Theorem 3.2}. For brevity denote $A:=A\left( \mathbf{x}%
_{0}\right) ,\tau :=\tau \left( \mathbf{x}_{0}\right) .$ Since by Theorem %
\ref{thm 1}, the function $\varphi \left( k\right) $ is analytic for $k\in 
\mathbb{R}$, then we assume below in this proof that the function $\varphi
\left( k\right) $ is given for all $k\in \mathbb{R}$. Using (\ref{4}), (\ref%
{5}) and (\ref{13}), we obtain for sufficiently large $k$%
\begin{equation}
\varphi \left( k\right) =A^{2}-2A\cos \left[ k\left( \tau -x_{3}\right) %
\right] +p\left( k\right) ,  \label{14}
\end{equation}%
where the real valued function $p\left( k\right) $ is such that 
\begin{equation}
p^{\left( j\right) }\left( k\right) =O\left( k^{-1}\right) ,\quad
k\rightarrow \infty ,j=0,1,2.  \label{15}
\end{equation}

It follows from (\ref{14}) and (\ref{15}) that $\displaystyle%
\lim_{k\rightarrow \infty }\varphi \left( k\right) $ exists if and only if $%
\tau =x_{3}.$ Hence, assume first that $\tau =x_{3}.$ To find the number $A$%
, we use the formula (4.16) of \cite[Theorem 1]{KlibanovRomanov:ip2016},%
\begin{equation}
A\left( \mathbf{x}_{0}\right) =\exp \left( -\frac{1}{2}\dint\limits_{\Gamma
\left( \mathbf{x}_{0},a\right) }\frac{1}{c\left( \mathbf{\xi }\right) }%
\Delta _{\xi }\tau \left( \mathbf{\xi }\right) d\sigma \right) .  \label{150}
\end{equation}%
By Lemma 2.1 $\Gamma \left( \mathbf{x}_{0},a\right) =L\left( \mathbf{x}%
_{0},a\right) .$ Furthermore, it follows from the assertion of that lemma
about $O\left( \mathbf{x}^{\prime }\right) $ that $\Delta _{\xi }\tau \left( 
\mathbf{\xi }\right) =0$ for all points $\mathbf{\xi }\in \Gamma \left( 
\mathbf{x}_{0},a\right) ,$ except of probably one point$.$ Hence, (\ref{150}%
) implies that $A\left( \mathbf{x}_{0}\right) =1.$

Consider now the case $\tau \neq x_{3}$. Denote 
\begin{equation}
\alpha \left( \mathbf{x}_{0}\right) =\alpha =\tau \left( \mathbf{x}%
_{0}\right) -x_{3}.  \label{151}
\end{equation}
Lemma 2.1 implies that 
\begin{equation}
\alpha >0.  \label{16}
\end{equation}%
By (\ref{14}) and (\ref{15}) 
\begin{equation}
\varphi ^{\prime }\left( k\right) =2A\alpha \sin \left( k\alpha \right)
+p^{\prime }\left( k\right) .  \label{17}
\end{equation}%
For sufficiently large $k$, consider the equation 
\begin{equation}
\varphi ^{\prime }\left( k\right) =0.  \label{18}
\end{equation}%
Consider the real valued function $q\left( k\right) ,$ 
\begin{equation*}
q\left( k\right) =-\frac{p^{\prime }\left( k\right) }{2A\alpha }.
\end{equation*}%
Then, using (\ref{15})-(\ref{17}), we obtain that equation (\ref{18}) is
equivalent with%
\begin{equation}
\sin \left( k\alpha \right) =q\left( k\right) ,  \label{19}
\end{equation}%
\begin{equation}
q\left( k\right) =O\left( \frac{1}{k}\right) ,q^{\prime }\left( k\right)
=O\left( \frac{1}{k}\right) ,k\rightarrow \infty .  \label{20}
\end{equation}%
Consider a sufficiently large integer $n>1.$ By (\ref{16}) we can choose a
sufficiently large $k>0$ such that 
\begin{equation}
k\alpha \in \left( \left( n-1\right) \pi ,n\pi +1\right)  \label{212121}
\end{equation}%
Since by (\ref{20}) $\left\vert q\left( k\right) \right\vert <1$ for
sufficiently large $k>0,$ then it follows from (\ref{19})-(\ref{212121}) that 
\begin{equation}
k\alpha =\left( -1\right) ^{n}\arcsin \left( q\left( k\right) \right) +n\pi .
\label{22}
\end{equation}

We show now that equation (\ref{22}) has unique solution 
\begin{equation*}
k\in \left( \frac{\left( n-1\right) \pi }{\alpha },\frac{n\pi +1}{\alpha }%
\right) =I_{n}.
\end{equation*}%
Indeed, consider the function $h_{n}\left( k\right) ,$ 
\begin{equation}
h_{n}\left( k\right) =k\alpha -\left( -1\right) ^{n}\arcsin \left( q\left(
k\right) \right) -n\pi ,k\in I_{n}.  \label{23}
\end{equation}%
By (\ref{20}) we can assume that 
\begin{equation}
\left\vert \arcsin \left( q\left( k\right) \right) \right\vert <\frac{1}{2}%
,k\in I_{n}.  \label{24}
\end{equation}%
Using (\ref{23}) and (\ref{24}), we obtain 
\begin{equation}
h_{n}\left( \frac{\left( n-1\right) \pi }{\alpha }\right) =-\pi -\left(
-1\right) ^{n}\arcsin \left( q\left( \frac{\left( n-1\right) \pi }{\alpha }%
\right) \right) <-2<0,  \label{25}
\end{equation}%
\begin{equation}
h_{n}\left( \frac{n\pi +1}{\alpha }\right) =1-\left( -1\right) ^{n}\arcsin
\left( q\left( \frac{n\pi +1}{\alpha }\right) \right) >\frac{1}{2}>0.
\label{26}
\end{equation}

It follows from (\ref{25}) and (\ref{26}) that the function $h_{n}\left(
k\right) $ has at least one zero inside of the interval $I_{n}.$ To show
that this zero is unique, consider the derivative $h_{n}^{\prime }\left(
k\right) ,$%
\begin{equation*}
h_{n}^{\prime }\left( k\right) =\alpha -\left( -1\right) ^{n}\frac{q^{\prime
}\left( k\right) }{\sqrt{1-q^{2}\left( k\right) }},k\in I_{n}.
\end{equation*}%
Since $n>1$ is sufficiently large and $k\in I_{n},$ we can assume by (\ref%
{20}) that 
\begin{equation*}
\left\vert \frac{q^{\prime }\left( k\right) }{\sqrt{1-q^{2}\left( k\right) }}%
\right\vert <\frac{\alpha }{2}.
\end{equation*}%
Hence, $h_{n}^{\prime }\left( k\right) >\alpha /2>0$ for $k\in I_{n}.$
Hence, the function $h_{n}\left( k\right) $ is monotonically increasing on
the interval $I_{n}.$ Hence, the above mentioned zero of the function $%
h_{n}\left( k\right) $ on the interval $I_{n}$ is unique. We denote this
zero $k_{n}.$

By (\ref{20})--(\ref{22})%
\begin{equation*}
k_{n}\alpha =n\pi +O\left( \frac{1}{n}\right) ,n\rightarrow \infty ,
\end{equation*}%
\begin{equation*}
k_{n+1}\alpha =\left( n+1\right) \pi +O\left( \frac{1}{n}\right)
,n\rightarrow \infty .
\end{equation*}%
Hence,%
\begin{equation}
\left( k_{n+1}-k_{n}\right) \alpha =\pi +O\left( \frac{1}{n}\right)
,n\rightarrow \infty .  \label{27}
\end{equation}%
In particular, it follows from (\ref{27}) that $k_{n+1}-k_{n}\geq \pi /2\neq
0$ for sufficiently large $n$. Thus, we obtain from (\ref{27})%
\begin{equation}
\alpha =\lim_{n\rightarrow \infty }\frac{\pi }{k_{n+1}-k_{n}}.  \label{28}
\end{equation}%
Since the function $\varphi ^{\prime }\left( k\right) $ is known, then all
its zeros are also known. Next, since equations (\ref{18}) and (\ref{19})
are equivalent, then (\ref{28}) uniquely defines the number $\tau =\alpha
+x_{3}$ as 
\begin{equation*}
\tau =x_{3}+\lim_{n\rightarrow \infty }\frac{\pi }{k_{n+1}-k_{n}}.
\end{equation*}%
Finally (\ref{15}) and (\ref{17}) imply that%
\begin{equation*}
A=\frac{1}{2\alpha }\lim_{n\rightarrow \infty }\varphi ^{\prime }\left[ 
\frac{1}{\alpha }\left( \frac{\pi }{2}+2n\pi \right) \right] .\text{ }
\end{equation*} The proof is complete. $\square$
%\end{proof}

\section{Numerical method for the approximate phase retrieval for $\mathbf{x}%
\in P_{\text{meas}}$}

\label{sec: phase retrieval}

In this section, we show how to recover functions $\tau (\mathbf{x})$ and $A(%
\mathbf{x})$ for $\mathbf{x}\in P_{\text{meas}}$ from the function $f\left( 
\mathbf{x},k\right) $ in (\ref{eqn f}). Everywhere in this section $\mathbf{x%
}\in P_{meas}.$ By Theorem \ref{thm 1} 
\begin{equation}
u(\mathbf{x},k)=A(\mathbf{x})\exp (\mathrm{i}k\tau (\mathbf{x}%
))(1+O(1/k)),k\rightarrow \infty ,  \label{eqn assymptotic}
\end{equation}%
where the function $u(\mathbf{x},k)$ is the solution of 
\eqref{eqn
Helmholtz}. In particular, this means that $\left\vert u(\mathbf{x}%
,k)\right\vert \approx A(\mathbf{x})$ for sufficiently large $k$. We show
below in this section how to approximate the functions $A(\mathbf{x})$ and $%
\tau (\mathbf{x})$ from the data $f(\mathbf{x},k)$, $\mathbf{x}\in P_{%
\mathrm{meas}},k\in \lbrack \underline{k},\overline{k}]$. Dropping the
remainder term $O(1/k)$ in \eqref{eqn assymptotic}, we deduce from %
\eqref{eqn total field} that 
\begin{equation}
u_{\mathrm{sc}}(\mathbf{x},k)=A(\mathbf{x})\exp (\mathrm{i}k\tau (\mathbf{x}%
))-\exp (\mathrm{i}kx_{3}).  \label{4.2}
\end{equation}%
Hence, for $\alpha (\mathbf{x})$ defined in (\ref{151}), the data is
approximated as 
\begin{equation}
f(\mathbf{x},k)=|u_{\mathrm{sc}}(\mathbf{x},k)|^{2}=A^{2}(\mathbf{x})+1-2A(%
\mathbf{x})\cos (k\alpha (\mathbf{x})).  \label{4.3}
\end{equation}

\begin{remark}[\protect\cite{KlibanovRomanov:SIAMam2016,
KlibanovRomanov:ip2016}]
Fix $\mathbf{x}$ in $P_{\mathrm{meas}}$. One can approximate $A(\mathbf{x})$
and $\alpha (\mathbf{x})$ by calculating the period of the function $%
k\mapsto f(\mathbf{x},k)$. For instance, we can find two consecutive local
minimizers (or maximizers) $\kappa _{1}$ and $\kappa _{2}$ of $f(\mathbf{x}%
,k)$ in $[\underline{k},\overline{k}]$. Thus, $\cos (\kappa _{1}\alpha (%
\mathbf{x}))=\cos (\kappa _{2}\alpha (\mathbf{x}))=1$ and 
\begin{equation}
\alpha (\mathbf{x})=\frac{2\pi }{|\kappa _{2}-\kappa _{1}|}.
\label{method Romanov}
\end{equation}%
This method can be used in theory. However, in some physical situations, the
interval $[\underline{k},\overline{k}]$  is not large enough for us to find
two local minima of $f(\mathbf{x},k)$. In addition, we might have errors
when finding local minimizers due to the noise added to the data, see Figure %
\ref{fig data} for an illustration. 
\begin{figure}[h]
\centering
\subfloat[Noiseless data.]{\includegraphics[width = .3\textwidth]{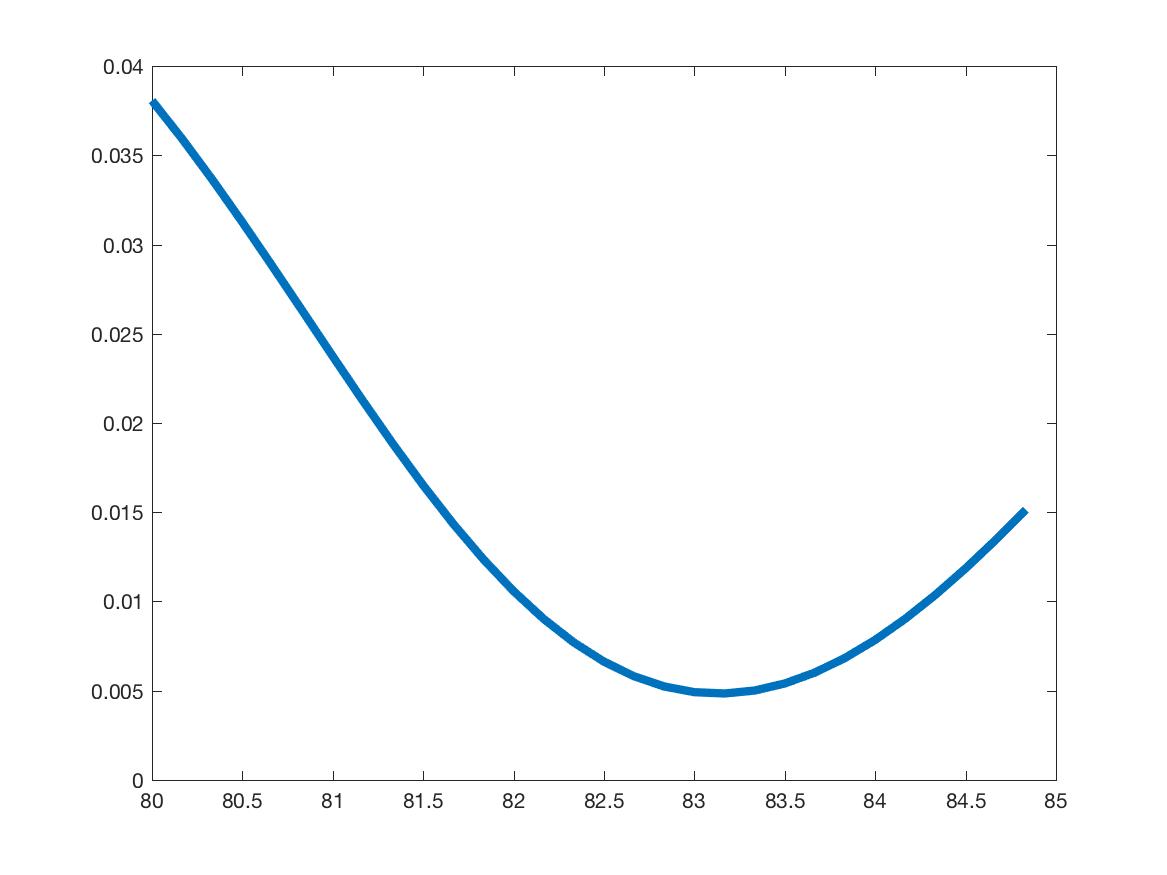}}
\hfill 
\subfloat[Data with 10\% noise]{\includegraphics[width =
.3\textwidth]{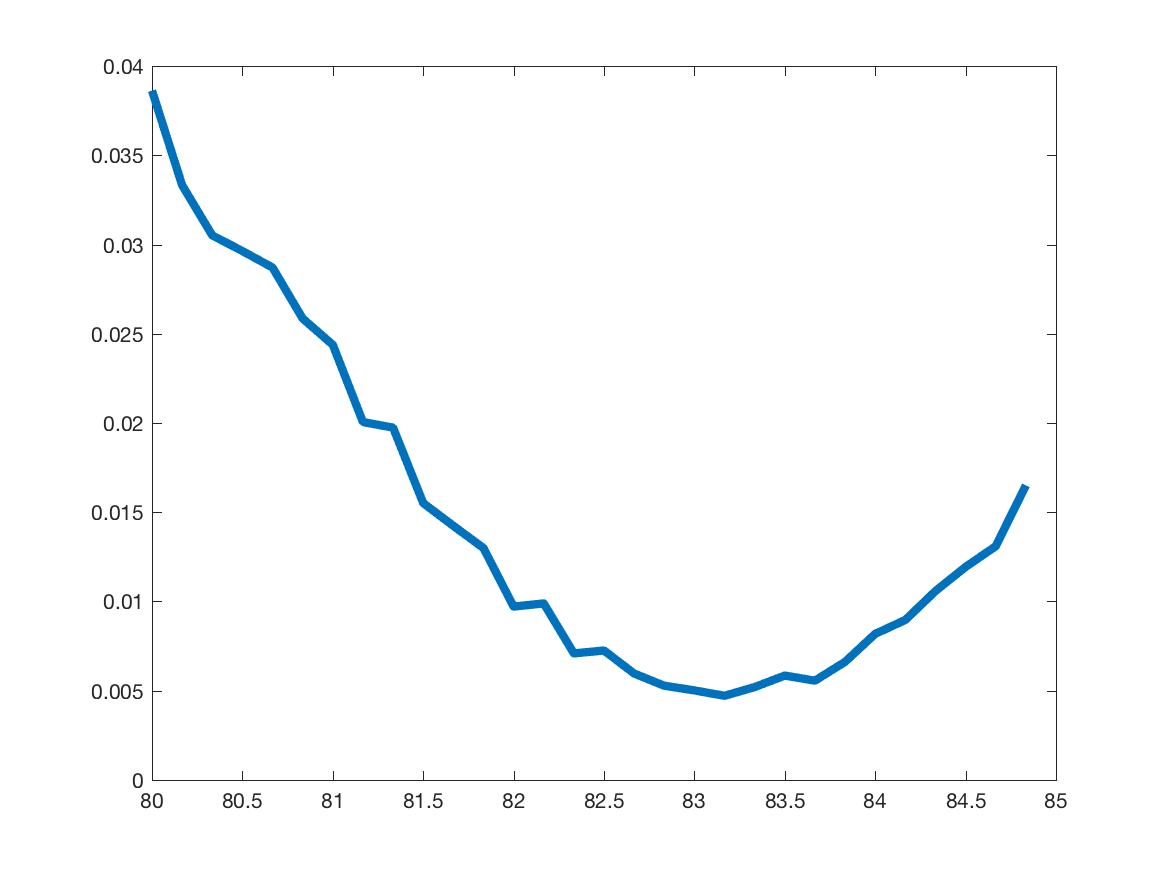}}\hfill 
\subfloat[Data with 15\%
noise]{\includegraphics[width = .3\textwidth]{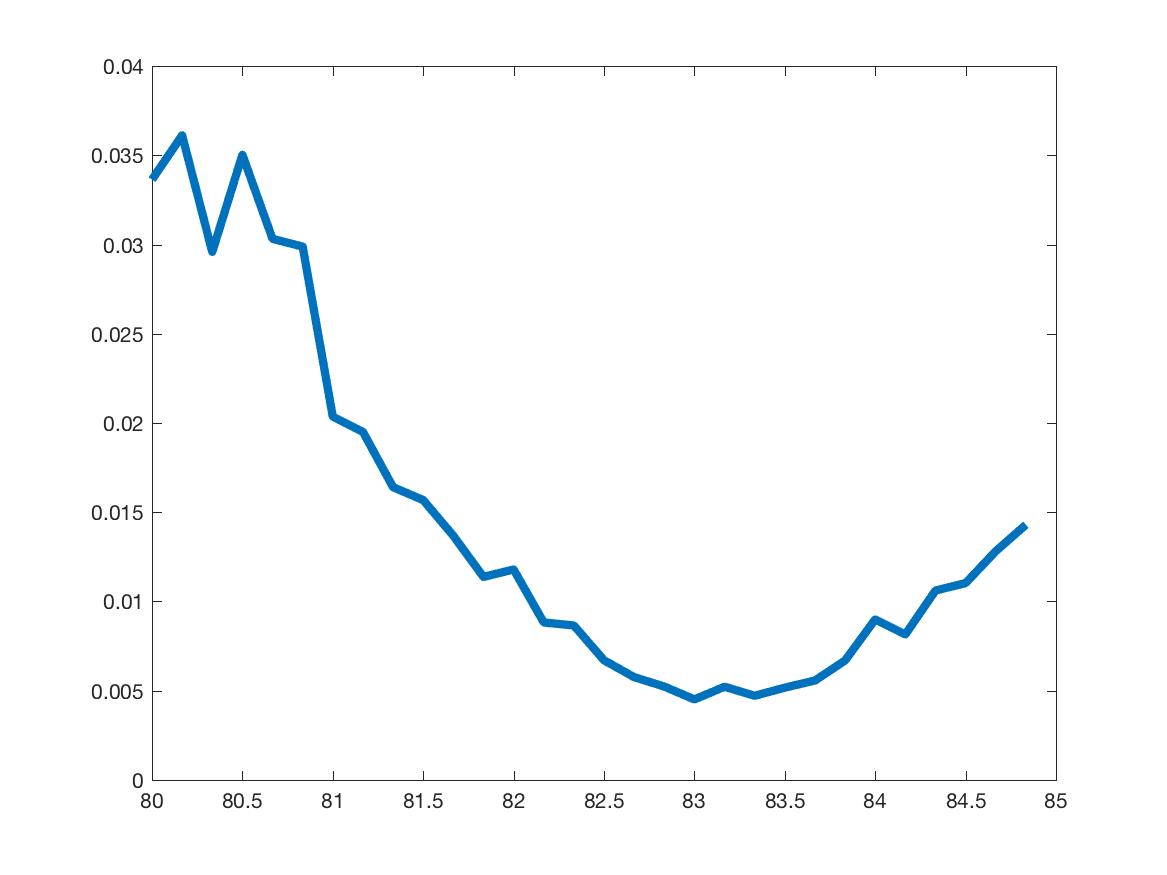}}
\caption{\textit{A typical example for the data $f(\mathbf{x},k)$ when $k$
varies in $[\protect\underline{k},\overline{k}] = [80, 85]$ where $[80, 85]$ is a realistic range of wave numbers, see Section \ref{sect:details}. In (a), the data attains
only one minimum value, which does not provide enough data to apply 
\eqref{method
Romanov}. The reconstruction is even more inconvenient when
the data, with noise, attains multiple extrema in the case (b) and (c).}}
\label{fig data}
\end{figure}
\end{remark}

%\begin{center}
%\textbf{Liem, Loc,}
%\end{center}
%
%\textbf{In label for Figure 2} \textquotedblleft in $\left[ \underline{k},%
%\overline{k}\right] ."$ \textbf{replace with} \textquotedblleft \emph{in }$%
%\left[ \underline{k},\overline{k}\right] =\left[ 80,85\right] ,$\emph{\
%where }$\left[ 80,85\right] $\emph{\ is a realistic range of wave numbers,
%see Section 6.3.}" \textcolor{red}{Done}

Modifying the phase reconstruction procedure of \cite%
{KlibanovLocKejia:apnum2016} where the incident wave is a point source
rather than the plane wave of our case, we propose the following
reconstruction process. Fix a point $\mathbf{x}\in P_{\mathrm{meas}}$ and
let $k\in \left[ \underline{k},\overline{k}\right] .$ Assuming that $\alpha (%
\mathbf{x})\not=0$, introduce $F_{1}(\mathbf{x},k)$ as 
\begin{align*}
F_{1}(\mathbf{x},k)& =\dint\limits_{\underline{k}}^{k}f(\mathbf{x},\kappa
)d\kappa =(A^{2}(\mathbf{x})+1)(k-\underline{k})-\frac{2A(\mathbf{x})}{%
\alpha (\mathbf{x})}(\sin (k\alpha (\mathbf{x}))-\sin (\underline{k}\alpha (%
\mathbf{x}))) \\
& =(A^{2}(\mathbf{x})+1)(k-\underline{k})+\frac{2A(\mathbf{x})}{\alpha (%
\mathbf{x})}\sin (\underline{k}\alpha (\mathbf{x}))-\frac{2A(\mathbf{x})}{%
\alpha (\mathbf{x})}\sin (k\alpha (\mathbf{x})).
\end{align*}%
Next, we define 
\begin{multline}
F_{2}(\mathbf{x},k)=\dint\limits_{\underline{k}}^{k}F_{1}(\mathbf{x},\kappa
)d\kappa =\frac{(A^{2}(\mathbf{x})+1)}{2}(k-\underline{k})^{2} \\
+\frac{2A(\mathbf{x})\sin (\underline{k}\alpha (\mathbf{x}))}{\alpha (%
\mathbf{x})}(k-\underline{k})+\frac{2A(\mathbf{x})}{\alpha ^{2}(x)}\cos
(k\alpha (x))-\frac{2A(\mathbf{x})}{\alpha ^{2}(x)}\cos (\underline{k}\alpha
(x)).  \label{4.5}
\end{multline}%
Combining \eqref{4.3} and \eqref{4.5} gives 
\begin{multline}
\alpha ^{2}(\mathbf{x})F_{2}(\mathbf{x},k)=\frac{\alpha ^{2}(\mathbf{x}%
)(A^{2}(\mathbf{x})+1)}{2}(k-\underline{k})^{2} \\
+2\alpha (\mathbf{x})A(\mathbf{x})\sin (\underline{k}\alpha (\mathbf{x}))(k-%
\underline{k})+A^{2}(\mathbf{x})+1-f(\mathbf{x},k)-2A(\mathbf{x})\cos (%
\underline{k}\alpha (x)).  \label{4.6}
\end{multline}%
Equation \eqref{4.6} can be rewritten as 
\begin{equation}
F_{2}(\mathbf{x},k)\xi _{1}(\mathbf{x})+(k-\underline{k})^{2}\xi _{2}(%
\mathbf{x})+(k-\underline{k})\xi _{3}(\mathbf{x})+\xi _{4}(\mathbf{x})=f(%
\mathbf{x},k)  \label{4.7}
\end{equation}%
for all $k\in \lbrack \underline{k},\overline{k}]$ where 
\begin{equation}
\begin{array}{rclrcl}
\xi _{1}(\mathbf{x}) & = & \alpha ^{2}(\mathbf{x}), & \xi _{2}(\mathbf{x}) & 
= & \displaystyle-\frac{\alpha ^{2}(\mathbf{x})(A^{2}(\mathbf{x})+1)}{2}, \\ 
\xi _{3}(\mathbf{x}) & = & -2\alpha (\mathbf{x})A(\mathbf{x})\sin (%
\underline{k}\alpha (\mathbf{x})), & \xi _{4}(\mathbf{x}) & = & A^{2}(%
\mathbf{x})-2A(\mathbf{x})\cos (\underline{k}\alpha (\mathbf{x}))+1.%
\end{array}
\label{4.70}
\end{equation}

Consider now the case $\alpha (\mathbf{x})=0.$ Then by (\ref{4.70}) $\xi
_{1}(\mathbf{x})=\xi _{2}(\mathbf{x})=\xi _{3}(\mathbf{x})=0$ and $\xi _{4}(%
\mathbf{x})=A^{2}(\mathbf{x})-2A\left( \mathbf{x}\right) +1.$ Hence, it
follows from (\ref{4.3}) that (\ref{4.7}) remains valid for the case $\alpha
(\mathbf{x})=0.$

For each $k\in \lbrack \underline{k},\overline{k}],$ equation \eqref{4.7} is
a linear equation with respect to the unknown vector $\mathbf{\xi (x)}%
=\left( {\xi }_{1}\mathbf{(x)},{\xi }_{2}(\mathbf{x}),{\xi }_{3}(\mathbf{x}),%
{\xi }_{4}\mathbf{(x)}\right) $ $.$ Consider the partition of the interval $%
\left[ \underline{k},\overline{k}\right] $ with the uniform step size $%
h=k_{j-1}-k_{j}\quad \mbox{for all }j\in \{1,\dots ,N\},$ 
\begin{equation}
k_{N}=\underline{k}<k_{N-1}<\dots <k_{1}<k_{0}=\overline{k}.  \label{203}
\end{equation}%
%
%
%
%
%
%Solving the over-determined linear system obtained by \eqref{4.7} with $k = k_i$, $i = \overline {0, N}$ provides $\alpha(\x) = \sqrt{|\xi_1(\x)|}$.
Then setting $k=k_{j}$ in (\ref{4.7}), we obtain a linear algebraic system
with respect to the vector $\mathbf{\xi .}$ This system is over-determined
since we have only 4 unknowns while $N+1$, the number of equations in the
system, is much greater than 4. We write this over-determined system as 
\begin{equation}
\mathcal{F}\mathbf{\xi }=\mathfrak{f},  \label{4.9}
\end{equation}%
where the $j^{\mathrm{th}}$row of the $\left( N+1\right) \times 4$ matrix $%
\mathcal{F}$ is given by $\left( F_{2}(\mathbf{x},k_{j}),(k_{j}-\underline{k}%
)^{2},(k_{j}-\underline{k}),1\right) $ and the $j^{\mathrm{th}}$ component
of the $N+1$ dimensional vector $\mathfrak{f}$ is $f(\mathbf{x},k_{j})$, $j=%
\overline{0,N}$. Then, we solve the following linear algebraic system 
\begin{equation}
(\mathcal{F}^{T}\mathcal{F}+\epsilon \mathrm{I}_{4})\xi =\mathcal{F}^{T}%
\mathfrak{f},  \label{4.10}
\end{equation}%
where $\mathcal{F}^{T}$ is $\mathcal{F}$ transpose and $\mathrm{I}_{4}$ is
the $4\times 4$ identity matrix. The positive small number $\epsilon $ plays
the role of regularization and its presence guarantees that \eqref{4.10} is
uniquely solvable. In our computations, the number $\epsilon $ is chosen by
a trial and error procedure (Section 6.3). After solving \eqref{4.10}, we
use \eqref{4.70} to set: 
\begin{equation}
\alpha (\mathbf{x})=\Re (\sqrt{\xi _{1}(\mathbf{x})}),\quad \tau (\mathbf{x}%
)=\alpha (\mathbf{x})+x_{3}.  \label{4.12}
\end{equation}%
Out of two possible values of $\sqrt{\xi _{1}(\mathbf{x})}$ we take the one
for which $\Re (\sqrt{\xi _{1}(\mathbf{x})})\geq 0.$ After obtaining $\alpha
(\mathbf{x})$, we compute $A(\mathbf{x})$ from~\eqref{4.3} as 
\begin{equation*}
A(\mathbf{x})=\left\vert \cos \left( k\alpha \left( \mathbf{x}\right)
\right) +\sqrt{\cos ^{2}\left( k\alpha \left( \mathbf{x}\right) \right) +f(%
\mathbf{x},k)-1}\right\vert .
\end{equation*}%
%
%
%
%We did not try to use the \textquotedblleft $+$" sign here since the \
%\textquotedblleft $-$" sign works well.
%
%\begin{center}
%\textbf{Loc }
%\end{center}
%2. How are we sure that \ $\Re \left[ \cos \left( k\alpha \left( \mathbf{x}%
%\right) \right) -\sqrt{\cos ^{2}\left( k\alpha \left( \mathbf{x}\right)
%\right) +f(\mathbf{x},k)-1}\right] >0$ Maybe we should use
%
%$\left\vert \cos \left( k\alpha \left( \mathbf{x}\right) \right) -\sqrt{\cos
%^{2}\left( k\alpha \left( \mathbf{x}\right) \right) +f(\mathbf{x},k)-1}%
%\right\vert ?$ Please comment in the text by BOLD\ FACED letters.
%
%\textbf{3. Do we use \textquotedblleft +" at the radical or
%\textquotedblleft }$-$\textbf{"?}
%
%{\bf From Loc: I leave this for Liem to answer.}
Let $\tau \left( \mathbf{x}\right) $ and $A(\mathbf{x})$ be two functions
reconstructed by the method of this section. Then following Theorem 3.1, we
obtain the following two approximate formulas for $\mathbf{x}%
=(x_{1},x_{2},x_{3})\in P_{\mathrm{meas}}:$ 
\begin{equation}
u\left( \mathbf{x},k\right) =A\left( \mathbf{x}\right) \exp \left( \text{i}%
k\tau \left( \mathbf{x}\right) \right) ,\quad u_{\mathrm{sc}}(\mathbf{x}%
,k)=A\left( \mathbf{x}\right) \exp \left( \text{i}k\tau \left( \mathbf{x}%
\right) \right) -\exp \left( \text{i}kx_{3}\right)  \label{4.14}
\end{equation}

\begin{figure}[h]
\centering
%\subfloat[Computed $A(\mathbf{x})$]{\includegraphics[width=0.35\textwidth]{Pic/A_1sphere}} 
%\subfloat[Computed
%$\tau(\x)$]{\includegraphics[width=0.35\textwidth]{Pic/tau_1sphere}}
%\par
\subfloat[The
real parts of the true (solid line) and
reconstructed (dashed line) scattered fields.]{\includegraphics[width=0.35\textwidth]{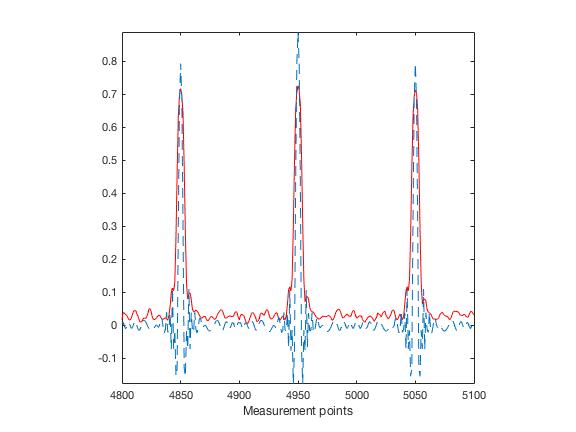}} \hspace{2cm}
\subfloat[The
imaginary parts of the true (solid line) and
reconstructed (dashed line) scattered fields.]{\includegraphics[width=0.35\textwidth]{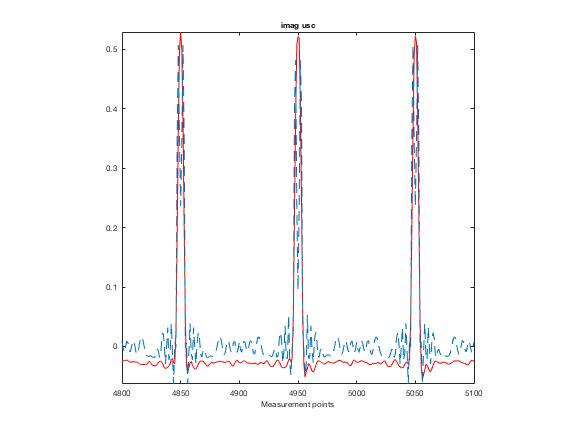}} 
\caption{\it  {An example of the
reconstructed functions} 
 $u_{sc}\left( \mathbf{x},k\right)$. {The function} $u_{sc}\left( \mathbf{x},k\right) $
{is computed by formula (\ref{4.14}) for} $k= 82.25$. {The data, with
5\% noise, for these computations correspond to Case 1 of Section \ref{sect:numerical}.} }
\label{fig reconstructed phase}
\end{figure}

To illustrate \eqref{4.14}, we arrange a uniform $100 \times 100$ grid points $\{\x_n\}_{n = 1}^{10,000}$ in $\Pm$ and show, in Figure \ref{fig reconstructed phase}, the true and reconstructed real and imaginary parts of $\usc(\x_n, k = 82.25)$ where $n \in \{4800, \dots, 5100\}$.

\begin{remark}
The reconstruction procedure described above is stable due to the stability
of the integration with respect to the noise. Indeed, our above analysis is
based on the integrals of the data.
\end{remark}

\begin{remark}[shifting the interval of wave numbers]
\label{rem shift} Formula {(\ref{4.14})} approximates $u_{\mathrm{sc}}(%
\mathbf{x},k)$ for $k\in \lbrack \underline{k},\overline{k}]$. However,
numerical solution of Helmholtz equation for large values of the wave number $%
k$ is very computationally expensive. Thus, in our numerical tests, we
extend {(\ref{4.14})} to another interval of wave numbers $k\in \lbrack 
\underline{k}^{\prime },\overline{k}^{\prime }]\subset \left( 0,\infty
\right) $ with smaller values of $k$, i.e. $\overline{k}^{\prime }<\overline{%
k}.$ To do this, we simply use in {(\ref{4.14}) values }$k\in \lbrack 
\underline{k}^{\prime },\overline{k}^{\prime }]${. }Thus, we first calculate 
$A\left( \mathbf{x}\right) ${\ and} $\tau \left( \mathbf{x}\right) $ using
values of $k\in \lbrack \underline{k},\overline{k}]$ from the original
interval and then use in {(\ref{4.14}) }$k\in \lbrack \underline{k}^{\prime
},\overline{k}^{\prime }]${\ . We assume everywhere below that this shift is
made and, to simplify notations, denote again }$\underline{k}^{\prime }:=%
\underline{k},\overline{k}^{\prime }:=\overline{k}.$
\end{remark}

In Sections 5 and 6, we briefly outline our globally convergent algorithm of 
\cite{Kliba2016}, which is playing an important role in our method to solve
the PCISP.

\section{The phased inverse scattering problem}

\label{sect: phased ICP}

We explain in Section 6.2 how to approximately obtain the boundary function $%
g\left( \mathbf{x},k\right) $ for $\mathbf{x}\in \partial \Omega ,k\in
\lbrack \underline{k},\overline{k}]$ in (\ref{eqn input for gca}) using {(%
\ref{4.14}). Hence, }our inverse scattering problem becomes now the phased
inverse scattering problem:

\begin{problem}[phased inverse scattering problem]
Given 
\begin{equation}
g(\mathbf{x},k)=u(\mathbf{x},k),\quad \mathbf{x}\in \partial \Omega ,k\in
\lbrack \underline{k},\overline{k}],  \label{202}
\end{equation}%
where $u(\mathbf{x},k)$ is the solution of \eqref{eqn Helmholtz}, determine
the function $c(\mathbf{x})$ for $\mathbf{x}\in \Omega .$ \label{problem
phased}
\end{problem}

The inverse problem~\eqref{202} has a broad range of applications and has
been widely studied. As to its uniqueness, we refer to item 2 in Remarks
2.2. We refer to~\cite%
{Ammar2004,Ammar2013,AmmariChowZou:sjap2016,BK,Colto1996,Gonch2017,Kliba2016,LiLiuWang:jcp2014, Li2015, Nov1,Nov2}
and references therein for various studies of numerical methods and
reconstruction procedures for solving this inverse problem under a variety
of assumptions on the measurement setup. The globally convergent algorithm
of \cite{Kliba2016} has been developed to solve the inverse problem~%
\eqref{202} with only a single measurement of multi-frequency scattering
data. In addition, we refer to \cite{KNT,convex1D} and references cited
therein for the second globally convergent numerical method for the single
measurement case, which is based on the construction of weighted globally
strictly convex Tikhonov-like functionals with Carleman weight functions in
them.

Below in this section we briefly describe the globally convergent numerical
method of \cite{Kliba2016}. We refer to \cite{Kliba2016} for details, which,
in particular, include the global convergence theorem 6.1.

\subsection{An integro-differential equation}

In this section, we assume that the function $u(\mathbf{x},k)$ never
vanishes. This assumption is true when $k$ is large due to 
\eqref{eqn
assymptotic}, see \cite{Kliba2016} for more details. Since the vector $%
\nabla u(\mathbf{x},k)/u(\mathbf{x},k)$ is curl free, we can follow a
procedure in \cite[Lemma 4.1]{Kliba2016} to find a smooth function $v(%
\mathbf{x},k)$ such that 
\begin{equation}
\exp (v(\mathbf{x},k))=u(\mathbf{x},k),\quad \nabla v(\mathbf{x},k)=\frac{%
\nabla u(\mathbf{x},k)}{u(\mathbf{x},k)}\quad \mathbf{x}\in \Omega ,k\in
\lbrack \underline{k},\overline{k}].  \label{log u}
\end{equation}%
The function $v$ can be understood as the natural logarithm of the function $%
u$. It satisfies 
\begin{equation}
\Delta v(\mathbf{x},k)+(\nabla v(\mathbf{x},k))^{2}=-k^{2}c(\mathbf{x}%
),\quad \mathbf{x}\in \Omega ,k\in \lbrack \underline{k},\overline{k}].
\label{eqn v}
\end{equation}%
Defining 
\begin{equation*}
q(\mathbf{x},k)=\frac{\partial v(\mathbf{x},k)}{\partial k},\quad \mathbf{x}%
\in \Omega ,k\in \lbrack \underline{k},\overline{k}]
\end{equation*}%
and differentiating \eqref{eqn v} with respect to $k$, we obtain that the
function $q(\mathbf{x},k)$ satisfies 
\begin{multline}
\frac{k}{2}\Delta q(\mathbf{x},k)+k\nabla q(\mathbf{x},k)\cdot \left(
-\dint\limits_{k}^{\overline{k}}\nabla q(\mathbf{x},s)s+\nabla V(\mathbf{x}%
)\right)
\\
=-\dint\limits_{k}^{\overline{k}}\Delta q(\mathbf{x},s)s+\Delta V(\mathbf{x}%
)+\left( -\dint\limits_{k}^{\overline{k}}\nabla q(\mathbf{x},s)s+\nabla V(%
\mathbf{x})\right) ^{2},\mathbf{x}\in \Omega ,k\in \lbrack \underline{k},%
\overline{k}]  \label{eqn ide}
\end{multline}%
and that, due to \eqref{log u}, $q(\mathbf{x},k)$ satisfies the following
Dirichlet boundary condition 
\begin{equation}
q(\mathbf{x},k)=\frac{\partial _{k}u(\mathbf{x},k)}{u(\mathbf{x},k)}\quad 
\mathbf{x}\in \partial \Omega .  \label{q boundary}
\end{equation}%
The function $V(\mathbf{x})$ is named the \textit{tail function}.

\subsubsection{The initial approximation $V_{0}( \mathbf{x}) $ for the tail
function}

\label{sec: tail}

Solving Problem \ref{problem phased} is somewhat equivalent to finding the
function $q(\mathbf{x},k)$, $\mathbf{x}\in \Omega $. Therefore, solving %
\eqref{eqn ide}--\eqref{q boundary} is crucial. However, the vector function 
$\nabla V(\mathbf{x})$ is involved in equation (\ref{eqn ide}), where the
tail function $V(\mathbf{x})$ is still unknown. In this subsection,
following \cite{Kliba2016}, we provide an initial approximation for $\nabla
V(\mathbf{x})$ and denote this approximation $\nabla V_{0}(\mathbf{x})$.
Thus, $\nabla V_{0}(\mathbf{x})$ is an important ingredient of the global
convergence theorem of \cite{Kliba2016}. We note that computing $\nabla
V_{0}(\mathbf{x})$ does not require any a priori knowledge of a good initial
guess for the true solution of Problem \ref{problem phased}. This is unlike
conventional locally convergent numerical methods.

Recall that the tail function is defined as $V(\mathbf{x})=v(\mathbf{x},%
\overline{k})$. Assuming that numbers $\underline{k}$ and $\overline{k}$ are
sufficiently large, dropping the term $O(1/k)$ in \eqref{eqn assymptotic}
and using (\ref{log u}), we obtain 
\begin{equation*}
V(\mathbf{x})\approx \ln A(\mathbf{x})+i\overline{k}\tau (\mathbf{x})=i%
\overline{k}\tau (\mathbf{x})\left( 1+\frac{\ln A(\mathbf{x})}{i\overline{k}%
\tau (\mathbf{x})}\right) \approx i\overline{k}\tau (\mathbf{x}),\quad 
\mathbf{x}\in \Omega .
\end{equation*}%
Therefore, the function $q(\mathbf{x},\overline{k})$ can be approximated as 
\begin{equation}
q(\mathbf{x},\overline{k})=\partial _{k}v(\mathbf{x},k)\mid _{k=\overline{k}%
}\approx i\tau (\mathbf{x})\approx \frac{V(\mathbf{x})}{\overline{k}},\quad 
\mathbf{x}\in \Omega .  \label{21}
\end{equation}%
Substituting \eqref{21} into \eqref{eqn ide} and setting in \eqref{eqn ide} $%
k=\overline{k}$, we obtain 
\begin{equation*}
\frac{1}{2}\Delta V(\mathbf{x})+(\nabla V(\mathbf{x}))^{2}=\Delta V(\mathbf{x%
})+\left( \nabla V(\mathbf{x})\right) ^{2},\quad \mathbf{x}\in \Omega ,
\end{equation*}%
which yields 
\begin{equation*}
\Delta V(\mathbf{x})=0,\quad \mathbf{x}\in \Omega .
\end{equation*}%
Note that only $\nabla V(\mathbf{x})$ and $\Delta V(\mathbf{x})=\mathrm{div}%
(\nabla V(\mathbf{x}))$ are involved in equation (\ref{eqn ide}). Hence,
rather that instead finding $V(\mathbf{x}),$ we compute directly the vector
function $\nabla V(\mathbf{x})$ in our numerical implementation. To do this,
we solve the following problem 
\begin{equation}
\left\{ 
\begin{array}{rcll}
\Delta (\nabla V(\mathbf{x})) & = & 0 & \mbox{in }\Omega , \\ 
\nabla V(\mathbf{x}) & = & R(\mathbf{x},\overline{k}) & \mbox{on }\partial
\Omega ,%
\end{array}%
\right.  \label{tail}
\end{equation}%
where $R(\mathbf{x},\overline{k})$ is a certain vector function, which is
known approximately. We refer to~\cite[Section 7.4]{Kliba2016} for all the
details about the approximation of $R(\mathbf{x},\overline{k})$ on the
entire boundary $\partial \Omega $ using, in particular, (\ref{202}) and (%
\ref{eqn input for gca}). We consider the solution of problem~\eqref{tail}
as the first approximation $\nabla V_{0}$ of the vector function $\nabla V$.
We mention once again that the globally convergent numerical method outlined
in this Section 5, including the approximation of \cite[Section 7.4]%
{Kliba2016} for the vector function $R(\mathbf{x},\overline{k}),$ has worked
quite well for the microwave experimental data, see references in Section 1.

\subsubsection{The globally convergent algorithm}

\label{sect:Algorithm1}

For $N\in \mathbb{N}$, consider the uniform partition 
\begin{equation}
k_{N}=\underline{k}<k_{N-1}<\dots <k_{1}<k_{0}=\overline{k}  \label{100}
\end{equation}%
of the interval $[\underline{k},\overline{k}]$ with the step size $%
h=k_{i-1}-k_{i}$, $1\leq i\leq N.$ Although this partition is different from
the one in (\ref{203}), we keep the same notation here for brevity. For each 
$n\in \{1,\dots ,N\},$ denote 
\begin{equation}
q_{n}(\mathbf{x})=q(\mathbf{x},k_{n}),\quad u_{n}(\mathbf{x})=u(\mathbf{x}%
,k_{n}),\quad \mathbf{x}\in \Omega .  \label{101}
\end{equation}%
Recall that the first approximation $\nabla V_{0}$ for the gradient $\nabla
V $ of the tail function is constructed in Section \ref{sec: tail}. We
assume, inductively, that $\nabla V_{n-1}$ is known, which implies that $%
\Delta V_{n-1}=\mathrm{div}\left( \nabla V_{n-1}\right) $ is known as well,
where $n\in \{1,\dots ,N\}$. By (\ref{100}) and (\ref{101}) the discrete,
with respect to $k$, form of equation \eqref{eqn
ide} is 
\begin{multline}
k_{n}\Delta q_{n}(\mathbf{x})-2k_{n}\nabla q_{n}(\mathbf{x})\cdot \nabla
Q_{n-1}(\mathbf{x})+2k_{n}\nabla q_{n}\cdot \nabla V_{n-1}(\mathbf{x}) \\
=-2\Delta Q_{n-1}(\mathbf{x})+2\Delta V_{n-1}(\mathbf{x})+2\left( -\nabla
Q_{n-1}(\mathbf{x})+\nabla V_{n-1}(\mathbf{x})\right) ^{2},\quad \mathbf{x}%
\in \Omega ,  \label{30}
\end{multline}%
where 
\begin{equation}
Q_{n-1}(\mathbf{x})=h\sum_{i=0}^{n-1}q_{n}(\mathbf{x}),\quad \mathbf{x}\in
\Omega .  \label{31}
\end{equation}%
Here, we approximate the integral $\displaystyle\int_{k}^{\overline{k}}q(%
\mathbf{x},s)ds$ by $Q_{n-1}(\mathbf{x})$ instead of $Q_{n}(\mathbf{x})$ to
remove the nonlinearity of \eqref{eqn ide}. The resulting error is $O(h),$
as $h\rightarrow 0.$ The boundary condition for the function $q_{n}(\mathbf{x%
})$ is 
\begin{equation}
q_{n}(\mathbf{x})=\frac{g(\mathbf{x},k_{n})-g(\mathbf{x},k_{n+1})}{hg(%
\mathbf{x},k_{n})},\quad \mathbf{x}\in \partial \Omega .  \label{32}
\end{equation}

\begin{remark}
Thus, (\ref{30})--(\ref{32}) is the Dirichlet boundary value problem for an
elliptic equation (\ref{30}). In \cite{Kliba2016}, for the theoretical
purpose, we make one more approximation for equation \eqref{30} via
replacing the term $2k_{n}\nabla q_{n}\cdot \nabla V_{n-1}(\mathbf{x})$ with
the term $2k_{n}\nabla q_{n-1}\cdot \nabla V_{n-1}(\mathbf{x}).$ The
resulting error is still $O(h)$ as $h\rightarrow 0.$ In this paper, although
we skip this approximation and use \eqref{30} to calculate the function $%
q_{n}$, the numerical results are still accurate, see Section \ref{sect: Num}%
.
\end{remark}

\begin{remark}
Although the derivative of the data with respect to $k$ is calculated in (%
\ref{32}) via the finite difference, we have not observed any instabilities
in our computations, probably because the step size $h$ was not exceedingly
small. The same is true for all above cited publications about the globally
convergent numerical methods of this group.
\end{remark}

In short, the algorithm is based on the following iterative process: (1)
given $\nabla V_{n-1}$, $\Delta V_{n-1},$ solve the Dirichlet boundary value
problem (\ref{30})--(\ref{32}); (2) use (\ref{eqn v}) at $k:=\underline{k}$
to calculate the function $c_{n}(\mathbf{x})$ via $q_{n}(\mathbf{x}),\nabla
V_{n-1}(\mathbf{x})$ and $\Delta V_{n-1}(\mathbf{x})$; (3) update $\nabla
V_{n}(\mathbf{x})$ and $\Delta V_{n}(\mathbf{x})$ by solving the
Lippmann-Schwinger equation (\ref{eqn LS}) with $\beta (\mathbf{x}):=c_{n}(%
\mathbf{x})-1$. To increase the stability of this iterative process, we
arrange internal iterations inside steps (1)-(3). The whole algorithm is
summarized as follows, see \cite{Kliba2016} for more details: 
\begin{algorithm}[H]
\caption{Globally convergent algorithm}
\label{alg:globalconv}
\begin{algorithmic}[1]
\State{Given $\nabla V_{0}$, set $q_{0}:=0$}
\For{$n = 1, 2, \dots, N$}
\State{Set $q_{n, 0} := q_{n - 1}$ and $\nabla V_{n, 0} := \nabla V_{n - 1}$}
\For{$i=1, 2, \dots, I_N$}
\State{Find $q_{n,i}$ by solving the elliptic boundary value problem~\eqref{30}--\eqref{32}}.
\State{Update 
$\nabla v_{n,i} :=- ( h\nabla q_{n,i}
+ h\sum_{j=0}^{n-1}\nabla q_{j}) +\nabla V_{n,i-1}$ in $\Omega$}.
\State{Update $c_{n,i}$ via~\eqref{eqn v}}.
\State{Find $u_{n,i}(\x,\overline{k})$ by solving the Lippmann-Schwinger  equation~\eqref{eqn LS} in $\Omega$ with $\beta(\mathbf{x}) :=c_{n,i}( \mathbf{x}) -1.$}
\State{Update  $\nabla V_{n,i}(\x) := \nabla u_{n,i}(\x,\overline{k})/ u_{n,i}(\x,\overline{k})$}.
\EndFor
\State{Update $q_{n} :=q_{n,I_N}$, $c_n :=c_{n,I_N}$ and $\nabla V_n := \nabla V_{n,I_N}$}.
\EndFor
\State{Choose $c$ by \emph{the-criterion-of-choice}}
\end{algorithmic}
\end{algorithm}

%\begin{center}
%\textbf{Loc,Liem, }
%\end{center}
%
%\textbf{In Algorithm 1:}
%
%\textbf{1. \textquotedblleft the boundary value problem" replace with
%\textquotedblleft the elliptic boundary value problem"}
%
%\textbf{2. I think that after Step 11 we should have command going to Step
%3? Also, we should say that the stopping criterion is not only for outer
%iterations but for inner iterations as well. True? I am NOT an expert in
%describing algorithms in this format. }

\begin{remark}
The stopping rule for the iterative loops in Algorithm~\ref{alg:globalconv} is presented 
in~\cite{Kliba2016, AlekKlibanovLocLiemThanh:anm2017, nguyen2017:iip2017,LiemKlibanovLocAlekFiddyHui:jcp2017}.
\end{remark}

\section{Numerical studies}

\label{sect: Num}

\subsection{Summary of our method for solving the PCISP}

In this section, we summarize the whole procedure of the reconstruction of
the coefficient $c(\mathbf{x})$ from the knowledge of $|u_{\mathrm{sc}}(%
\mathbf{x},k)|$, $\mathbf{x}\in P_{\mathrm{meas}}$ and $k\in \lbrack 
\underline{k},\overline{k}]$ as follows: 
\begin{algorithm}[H]
\caption{Globally convergent algorithm for the phaseless inverse scattering problem}
\begin{algorithmic}[1]
\State{For each point $\x$ in $ P_{\mathrm{meas}}$ find numbers $\tau \left( \mathbf{x}\right) $ and $A\left( \mathbf{x}\right) $ as described in section \ref{sec: phase retrieval}. Next, approximate the
function $u_{\text{sc}}\left( \mathbf{x},k\right) $ via (\ref{4.14}) for $\mathbf{x}\in P_{\mathrm{meas}}, k\in \left[ \underline{k},\overline{k}\right]
,$ where $\left[ \underline{k},\overline{k}\right] $ is the
``shifted" interval of wave numbers as in Remark \ref{rem
shift}.}
\State{\label{step propagation}Propagate the reconstructed $\usc(\x, k)$ from $\Pm$ to the plane  $\Pp$, see Section \ref{sec: data propagation}. The plane $\Pp$ is closer to the targets than $\Pm$.}
\State{Consider the square  $\Gamma \subset
\partial \Omega $ in (\ref{204204}), which is a part of the propagated plane $%
P_{\text{prop}}$, where $\Omega $ is the domain of our interest containing
all targets. Assign the data on $\partial \Omega $ as in (\ref{eqn input for
gca}).}
\State{Having approximated the
function }$u(\mathbf{x},k)$ on $\partial \Omega $ {as in (%
\ref{eqn input for gca}), find }$c\left( \mathbf{x}\right) ${ by
Algorithm \ref{alg:globalconv}.} 
\end{algorithmic}
\label{alg phaseless}
\end{algorithm}

%\begin{center}
%\textbf{Liem, Loc,}
%\end{center}
%
%\textbf{In Algorithm 2:}
%
%\textbf{1. \textquotedblleft For each }$\mathbf{x}$\textbf{\ in }$P_{\text{%
%meas}}"$\textbf{\ replace with \textquotedblleft For each point }$\mathbf{x}%
%\in P_{\text{meas}}"$ \textcolor{red}{Done.}
%
%\textbf{2. \textquotedblleft wavenumbers" replace with \textquotedblleft
%wave numbers" since we use this everywhere above} \textcolor{red}{Done.}
%
%\textbf{3. In Step 2 \textquotedblleft the }$P_{\text{prop}}"$\textbf{\
%replace with \textquotedblleft the plane }$P_{\text{prop}}"$ \textcolor{red}{Done.}
%
%\textbf{4. In Step 3 Replace \textquotedblleft (\ref{4.14})" with
%\textquotedblleft (6.2)"} \textcolor{red}{Done.}
%
%\textbf{5. Step 4 should be \textquotedblleft Having approximated the
%function }$u(\mathbf{x},k)$\textbf{\ on }$\partial \Omega $\textbf{\ as in (%
%\ref{eqn input for gca}), find }$c\left( \mathbf{x}\right) $\textbf{\ by
%Algorithm 1."}  \textcolor{red}{Done.}

\subsection{Data propagation and completion}

\label{sec: data propagation}

The data propagation is a procedure which enables us to \textquotedblleft
move" the data to a plane which is closer to the target than the original
measurement plane. Thus, we \textquotedblleft propagate" the reconstructed
function $u_{\mathrm{sc}}$ in (\ref{4.14}) from the square $P_{\mathrm{meas}%
}\subset P$ to $P_{\mathrm{prop}}$. Here $P_{\mathrm{meas}}$ is the square (%
\ref{measurement plane}) on the measurement plane $P=\left\{ \mathbf{x}%
:x_{3}=R\right\} $ and $P_{\mathrm{prop}}$ is a propagated plane which is
closer to the targets of our interest. This data propagation process has
been rigorously justified in~\cite{nguyen2017:iip2017}. By our experience in
the previous works~\cite{AlekKlibanovLocLiemThanh:anm2017,
nguyen2017:iip2017, LiemKlibanovLocAlekFiddyHui:jcp2017}, this process
enables one not only to propagate the scattered wave but also to
significantly decrease the amount of noise in the data. Moreover, unlike the
measured data, the propagated data focuses more at the $x_{1},x_{2}$
positions of the targets. We briefly outline the data propagation method
here.

Let the number $R^{\prime }\in \left( 0,R\right) .$ Assume that the domain $%
\Omega \subset \left\{ x_{3}\in \left( 0,R^{\prime }\right) \right\} $. So,
we want to propagate the function $u_{sc}\left( \mathbf{x},k\right) $ given
in (\ref{4.14}) for $k\in \lbrack \underline{k},\overline{k}]$ from $P_{%
\text{meas}}$ to the plane $P_{\text{prop}}=\left\{ x_{3}=R^{\prime
}\right\} .$ For any pair of real numbers $k_{x_{1}},k_{x_{2}}$, define 
\begin{equation*}
\widehat{u}_{\mathrm{sc}}(k_{x_{1}},k_{x_{2}},k)=\frac{1}{2\pi }%
\dint\limits_{\mathbb{R}^{2}}u_{\mathrm{sc}}(x_{1},x_{2},R)\exp (\mathrm{i}%
(k_{x_{1}}x_{1}+k_{x_{2}}x_{2}))dx_{1}dx_{2}.
\end{equation*}%
Here, we have extended $u_{\mathrm{sc}}(\mathbf{x},k)$ by zero for $\mathbf{x%
}=(x_{1},x_{2},R)\not\in P_{\mathrm{meas}}.$ Then it was proved in~\cite%
{nguyen2017:iip2017} that the propagated wave field $u_{\mathrm{sc}}(\mathbf{%
x},k)$ for $0<R^{\prime }<R$ is given by 
\begin{equation}
u_{\mathrm{sc}}(\mathbf{x},k)=\frac{1}{2\pi }\dint\limits_{%
\{k_{x_{1}}^{2}+k_{x_{2}}^{2}<k^{2}\}}\widehat{u}_{\mathrm{sc}%
}(k_{x_{1}},k_{x_{2}},k)\exp (-\mathrm{i}%
(k_{x_{1}}x_{1}+k_{x_{2}}x_{2}-k_{x_{3}}(R^{\prime
}-R)))dk_{x_{1}}dk_{x_{2}},  \label{5.1}
\end{equation}%
where $\mathbf{x}=(x_{1},x_{2},R^{\prime })$, $%
k_{x_{3}}=(k^{2}-k_{x_{1}}^{2}-k_{x_{2}}^{2})^{\frac{1}{2}}$ and $k\in \left[
\underline{k},\overline{k}\right] .$ We use the same square on the plane $P_{%
\text{prop}}$ as in $P_{\text{meas}}$ (see (\ref{measurement plane})) and we
do not count values of the function $u_{\mathrm{sc}}(\mathbf{x},k)$ in (\ref%
{5.1}) outside of this square. Thus, we denote that square on $P_{\text{prop}%
}$ as 
\begin{equation}
\Gamma =\left\{ \mathbf{x}:\left\vert x_{1}\right\vert <b,\left\vert
x_{2}\right\vert <b,x_{3}=R^{\prime }\right\} .  \label{204204}
\end{equation}%
We assume that $\Gamma \subset \partial \Omega .$

The data for our globally convergent algorithm are $u(\mathbf{x}%
,k)|_{\partial \Omega \times \lbrack \underline{k}^{\prime },\overline{k}%
^{\prime }]}$ \cite{Kliba2016}. Therefore, we need to complement the data on 
$\partial \Omega \setminus \Gamma $ as it was done in \cite{Kliba2016}. We
are doing so heuristically by simply setting $u_{\mathrm{sc}}(\mathbf{x}%
,k)=0 $ for $\mathbf{x}\in \partial \Omega \setminus \Gamma $. In other
words, the input $u(\mathbf{x},k)\mid _{\partial \Omega }:=g(\mathbf{x},k)$
for the globally convergent numerical method of \cite{Kliba2016} described
above is given by 
\begin{equation}
g(\mathbf{x},k)=\left\{ 
\begin{array}{ll}
u_{\mathrm{sc}}(\mathbf{x},k)+\exp (\mathrm{i}kx_{3}), & \mathbf{x}\in
\Gamma ,k\in \lbrack \underline{k},\overline{k}], \\ 
\exp (\mathrm{i}kx_{3}), & \mathbf{x}\in \partial \Omega \setminus \Gamma
,k\in \lbrack \underline{k},\overline{k}],%
\end{array}%
\right.  \label{eqn input for gca}
\end{equation}%
see Remark 4.3 for $[\underline{k},\overline{k}].$

\begin{remark}
{It was shown in subsections 7.6 and 7.7 of \cite{Kliba2016} that in the
case when the correct computationally simulated data are assigned on the
entire boundary }$\partial \Omega ,$ {\ the computational result is about
the same as for the case when the boundary data given on a part of the
boundary are complemented as in (\ref{eqn input for gca}). Also, it was
demonstrated in all our above cited works on experimental data that (\ref%
{eqn input for gca}) works well.}
\end{remark}

\subsection{Some details of numerical experiments}
\label{sect:details}

Our numerical studies are conducted for a realistic range of parameters
which we have extensively discussed with Professor Vasily Astratov from
Center for Optoelectronics and Optical Communications of the University of North
Carolina at Charlotte.

In at least one experimental arrangement one wants to image dielectric balls
whose diameters are about $5\mu m.$ These balls are called \textquotedblleft
microspheres". In our computations, the measurement plane is about $25\mu m$
away from the domain $\Omega $ where these microspheres are located. The
size of the measurement square $P_{\mathrm{meas}}$ is $100\mu m\times 100\mu
m.$ The wavelengths of light $\lambda \in \left[ 738,785\right] nm=\left[
0.738,0.785\right] \mu m.$

To make variables dimensionless, we consider the change of variables $%
\mathbf{x}^{\prime }=\mathbf{x/}10\mu m$ while leaving the same notations
for brevity. Then the dimensionless wave number $k=20\pi /\lambda ^{\prime },$
where $\lambda ^{\prime }$ is the dimensionless wavelength. Hence, we obtain
the interval for the dimensionless $k\in \left[ 80,85\right] =\left[ 20\pi
/0.785,20\pi /0.738\right] .$ We then \textquotedblleft shift" the interval
of wavelength to $[20.4,21]$ as in Remark \ref{rem shift} to make the whole
procedure less computationally expensive. More precisely, we consider the
following setup:

\begin{enumerate}
\item[(a)] The scattering balls are located near the $x_{1}x_{2}-$plane and
their diameter is 0.5.

\item[(b)] The domain $\Omega =(-2.5, 2.5) \times (-2.5, 2.5) \times (-4,
1). $

\item[(c)] The measurement square is $P_{\mathrm{meas}}=\{\mathbf{x}%
:|x_{1}|,|x_{2}|\leq 5,x_{3}=2.5\}.$

\item[(d)] The propagated plane in Step \ref{step propagation} in Algorithm %
\ref{alg phaseless} is \[P_{\mathrm{prop}}=\left\{\x = (x_1, x_2, x_3): |x_1|, |x_2| \leq 5, x_{3}=1\right\}. \]
\end{enumerate}

Then, we choose $\Gamma$ to be a subset of $\Pp$ as 
\begin{equation*}
\Gamma =\{\mathbf{x}=(x_{1},x_{2},x_{3}):|x_{1}|,|x_{2}|\leq 2.5,x_{3}=1\}.
\end{equation*} 
The main reason for working with $\Gamma$ instead of $P_{\mathrm{prop}}$
is that the data on $P_{\mathrm{prop}}\setminus \Gamma$ are small 
and do not contribute the inversion process. Furthermore, this choice leads to a smaller  computational
domain $\Omega$.
Thus, in notations of Section 6.2, $R=2$ and $R^{\prime }=1.$ We had $%
100\times 100$\ uniform grid $\left\{ \mathbf{x}_{n}\right\} $\ covering the
square $P_{\mathrm{meas}}$. Functions $A$ and $\tau $ were reconstructed at
these grid points. We have chosen the regularization parameter $\epsilon
=0.03$ in (\ref{4.10}) by trial and error. Thus, we have obtained numbers $%
A\left( \mathbf{x}_{n}\right) $ and $\tau \left( \mathbf{x}_{n}\right) .$%
These numbers were reconstructed from noisy data. For $k\in \left[ 80,85%
\right] $ the 5\% random noise was introduced as: 
\begin{equation*}
f_{\mathrm{noise}}(\mathbf{x},k)=f(\mathbf{x},k)+5\%\Vert f\Vert _{L^{2}(P_{%
\mathrm{meas}}\times \lbrack \underline{k},\overline{k}])}\mathrm{rand(}%
\mathbf{x}\mathrm{,}k\mathrm{)}/\Vert \mathrm{rand}\Vert _{L^{2}(P_{\mathrm{%
meas}}\times \lbrack \underline{k},\overline{k}])}
\end{equation*}%
for $\mathbf{x}\in P_{\mathrm{meas}},k\in \lbrack \underline{k},\overline{k}%
] $ where $\mathrm{rand(}\mathbf{x}\mathrm{\mathrm{,}}k\mathrm{)}$ is a
random number in $(0,1).$

\subsection{Numerical results}
\label{sect:numerical}

To make our spherical inclusions to be smoothly embedded in the background
medium, we consider the following construction. Let $B(\mathbf{x}_{0},r)$ be
the ball of the radius $r$ centered at the point $\mathbf{x}_{0}\in \mathbb{R%
}^{3}$. Define the function $\chi _{B(\mathbf{x}_{0},r)}(\mathbf{x})$, 
\begin{equation*}
\chi _{B(\mathbf{x}_{0},r)}(\mathbf{x})=\left\{ 
\begin{array}{ll}
\exp (1-r^{2}/(r^{2}-|\mathbf{x}-\mathbf{x}_{0}|^{2})) & \mathbf{x}\in B(%
\mathbf{x}_{0},r) \\ 
0 & \mathbf{x}\in \mathbb{R}^{3}\setminus \overline{B(\mathbf{x}_{0},r)}.%
\end{array}%
\right.
\end{equation*}%
Hence, the function $\chi _{B(\mathbf{x}_{0},r)}\in C^{\infty }(\mathbb{R}%
^{3})$ and its support is $B(\mathbf{x}_{0},r)$. We present below the
following three cases of the numerical reconstruction of the functions $c(%
\mathbf{x})$:

\begin{enumerate}
\item One inclusion: $c(\mathbf{x})=1+\chi _{B(\mathbf{x}_{0},r)}$ where $%
r=0.25$ and $\mathbf{x}_{0}=(0,0,0.25)$.

\item Two inclusions, which are symmetric with respect to the plane $%
\{x_{1}=0\}$: $c(\mathbf{x})=1+\chi _{B(\mathbf{x}^{\left( 1\right)
},r)}+\chi _{B(\mathbf{x}^{\left( 2\right) },r)}$ where $r=0.25$, $\mathbf{x}%
^{\left( 1\right) }=(-0.5,0,0.25)$ and $\mathbf{x}^{\left( 2\right)
}=(0.5,0,0.25)$.

\item Two inclusions which are non-symmetric with respect to the plane $%
\{x_{2}=0\}$ but symmetric with respect to the plane $\{x_{1}=0\}$: $c(%
\mathbf{x})=1+\chi _{B(\mathbf{x}^{\left( 3\right) },r)}+\chi _{B(\mathbf{x}%
^{\left( 4\right) },r)}$ where $r=0.25$, $\mathbf{x}^{\left( 3\right)
}=(0.5,0.5,0.25)$ and $\mathbf{x}^{\left( 4\right) }=(-0.5,-0.25,0.25)$.
\end{enumerate}

\begin{figure}[h!]
\centering
\subfloat[Exact profile (cross-sectional
view)]{\includegraphics[width=0.38\textwidth]{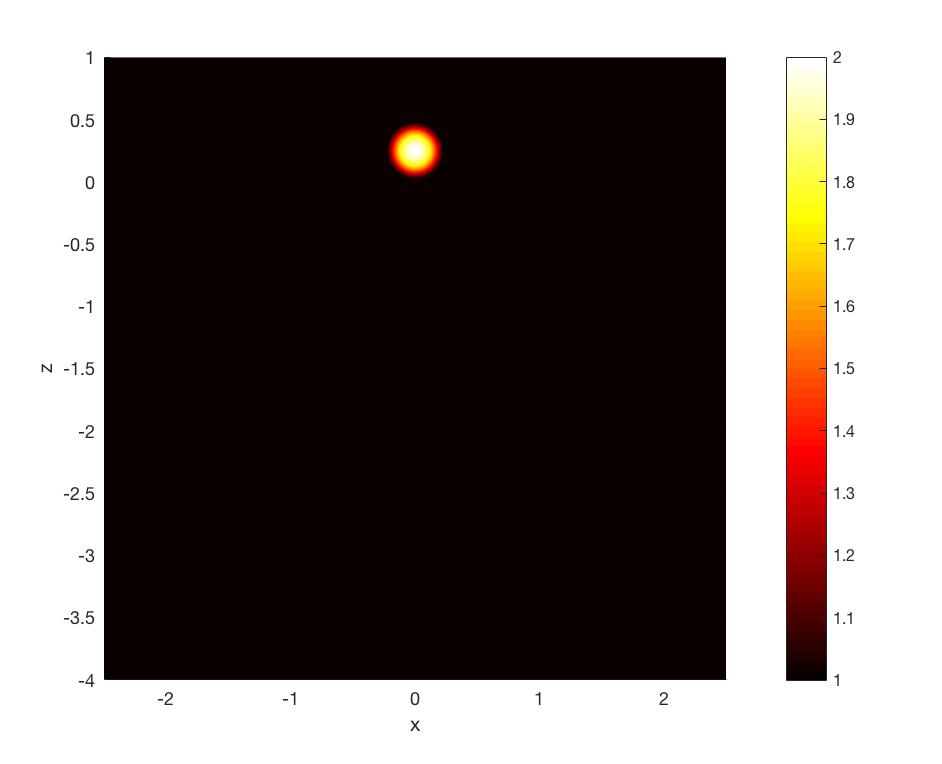}} 
\subfloat[Exact profile (3D
view)]{\includegraphics[width=0.4\textwidth]{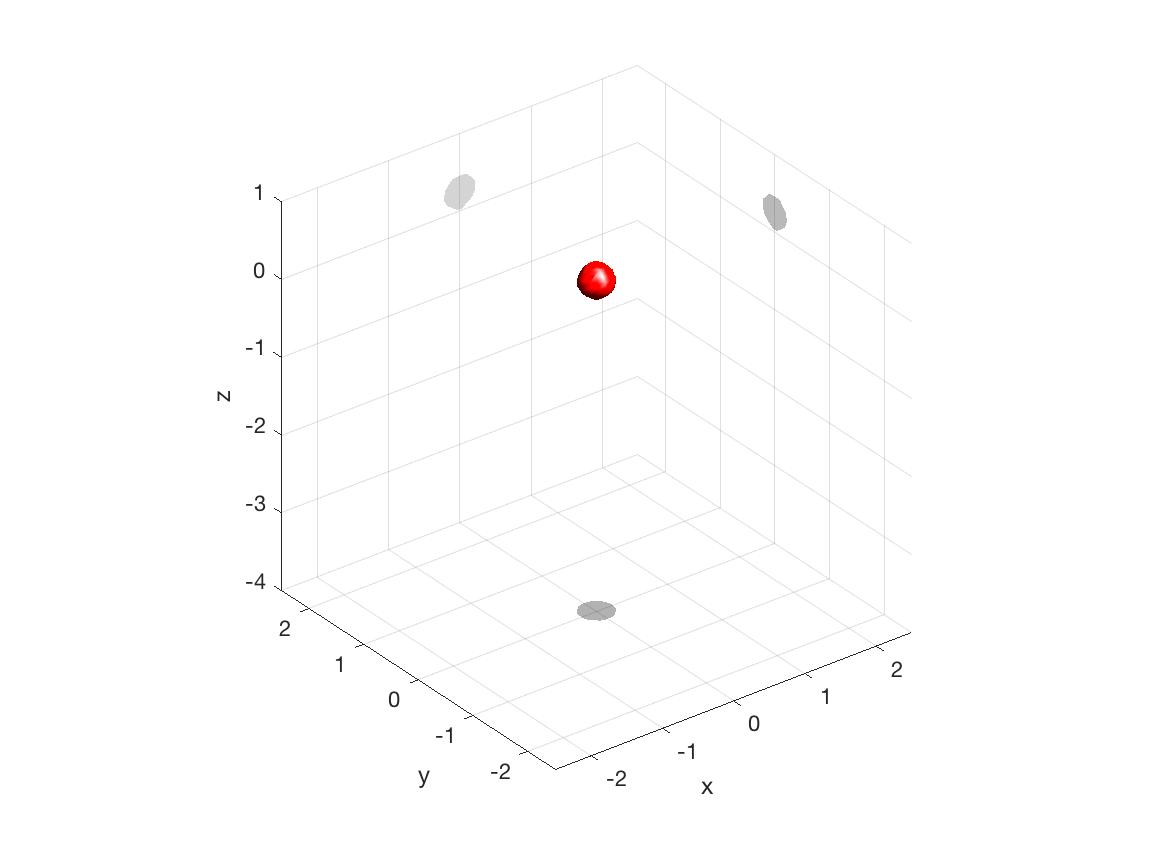}} \\
\subfloat[Reconstruction (cross-sectional
view)]{\includegraphics[width=0.4\textwidth]{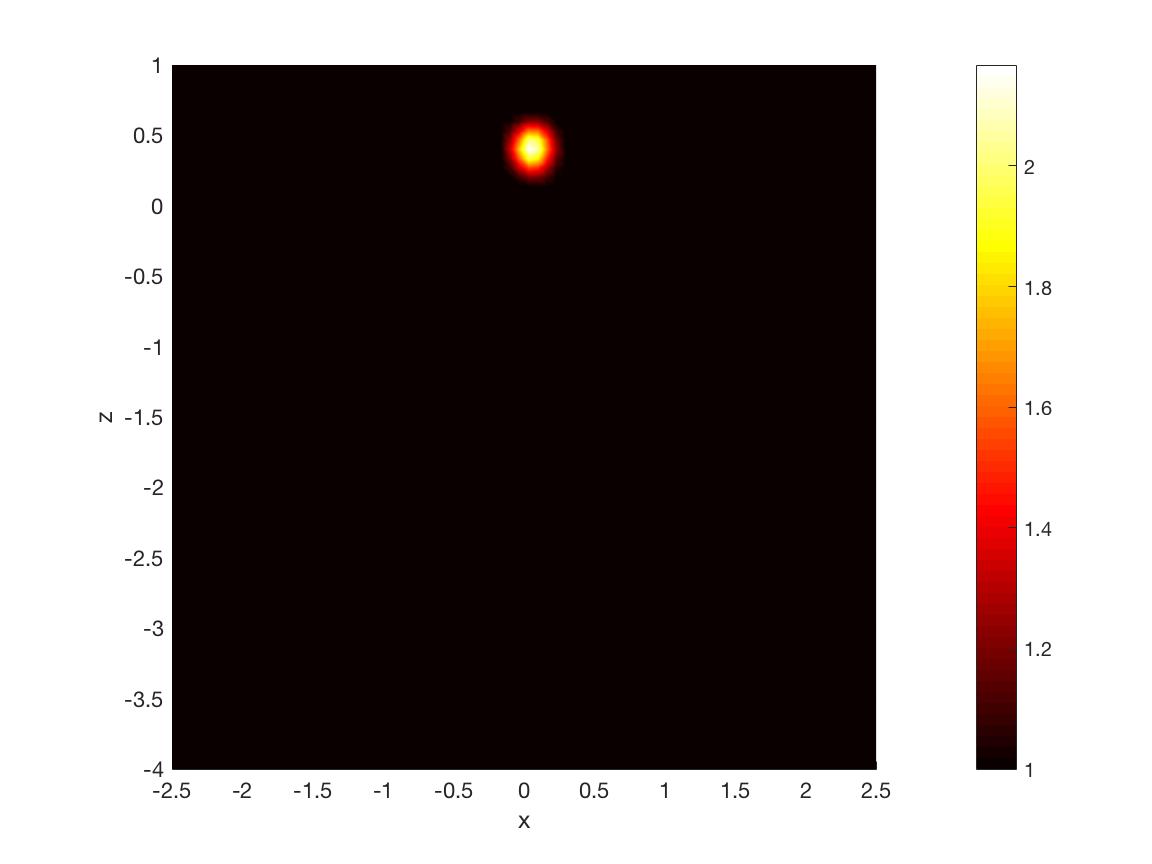}} 
\subfloat[Reconstruction (3D
view)]{\includegraphics[width=0.4\textwidth]{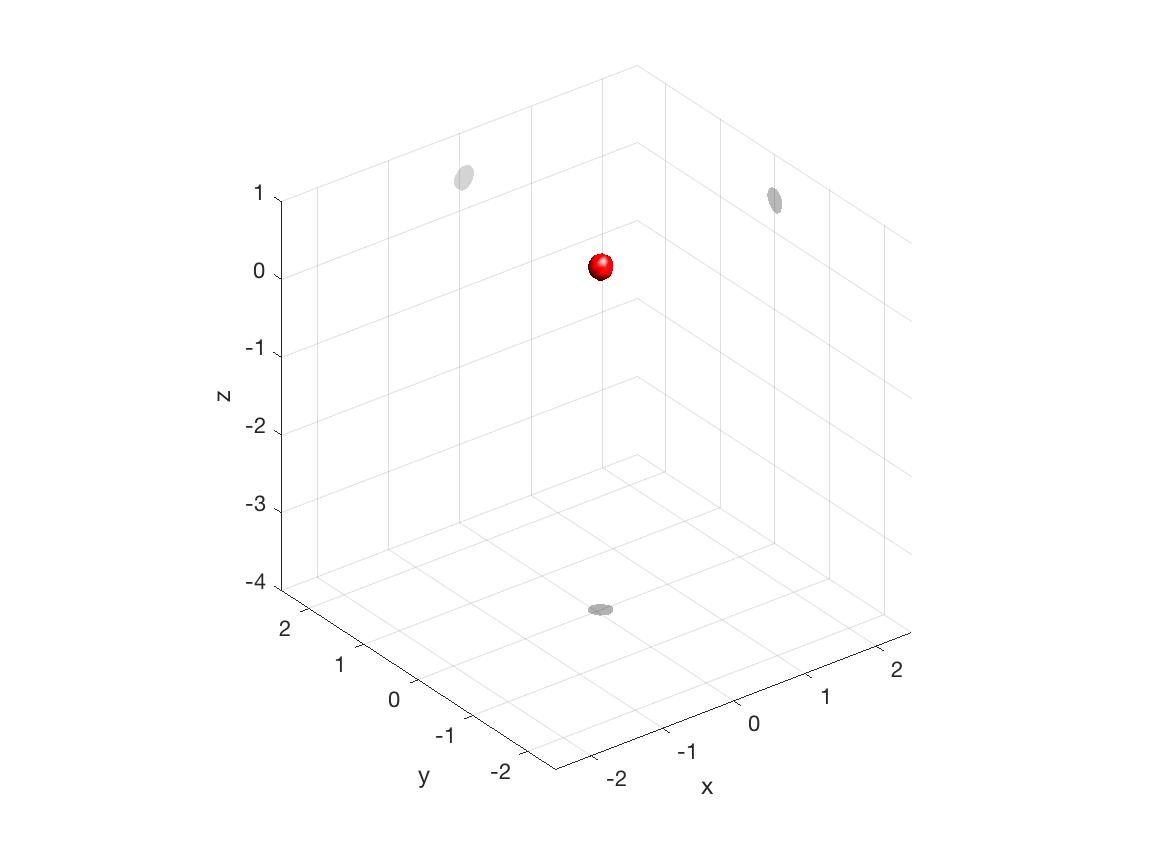}}
\caption{\textit{Reconstruction results for the case of one spherical inclusion: the above item 1. {The
maximal value of }$c^{\ast }\left( \mathbf{x}\right) ${\ in this inclusion is
2. The maximal value of }$c_{\text{comp}}\left( \mathbf{x}\right) ${\ in
this inclusion is 2.16. Hence, the error in computing this value is 8\%}}}
\label{fig:res_obj1}
\end{figure}

Here, an inclusion means a connected component of the support of the
function $c(\mathbf{x})-1.$ In labels for Figures \ref{fig:res_obj1}--\ref%
{fig:res_obj3} $c^{\ast }\left( \mathbf{x}\right) $ and $c_{\text{comp}%
}\left( \mathbf{x}\right) $ mean the exact and computed coefficients $c(%
\mathbf{x})$, respectively. We display in these figures:

\begin{enumerate}
\item[(a)] A 2D cross-sectional view of the true inclusions through their
center by a plane which is orthogonal to the $x_{1},x_{2}-$plane.

\item[(b)] A 3D view of the true inclusions by isosurfaces.

\item[(c)] A 2D cross-sectional view of the reconstructed inclusions on the
same plane as in (a).

\item[(d)] A 3D view of the reconstructed inclusions by isosurfaces.
\end{enumerate}

\begin{figure}[h!]
\centering
\subfloat[Exact profile (cross-sectional
view)]{\includegraphics[width=0.4\textwidth]{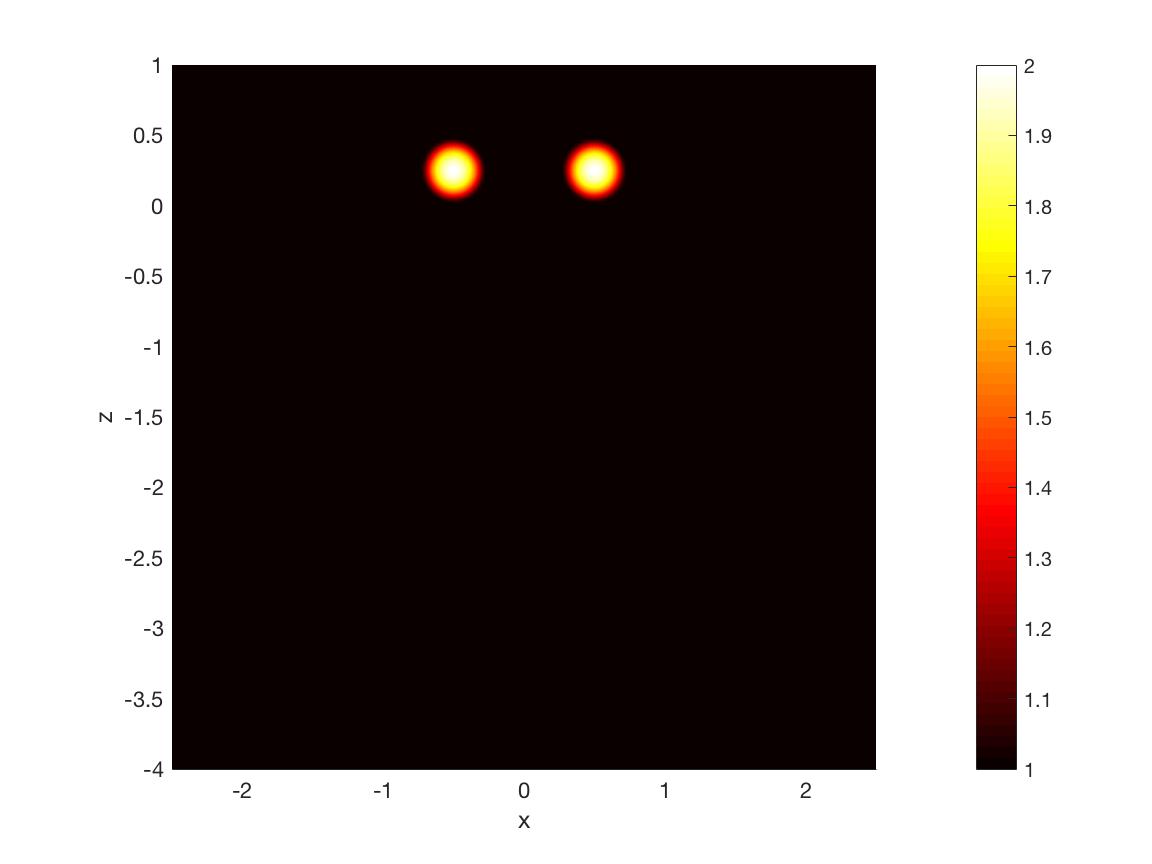}} 
\subfloat[Exact profile (3D
view)]{\includegraphics[width=0.4\textwidth]{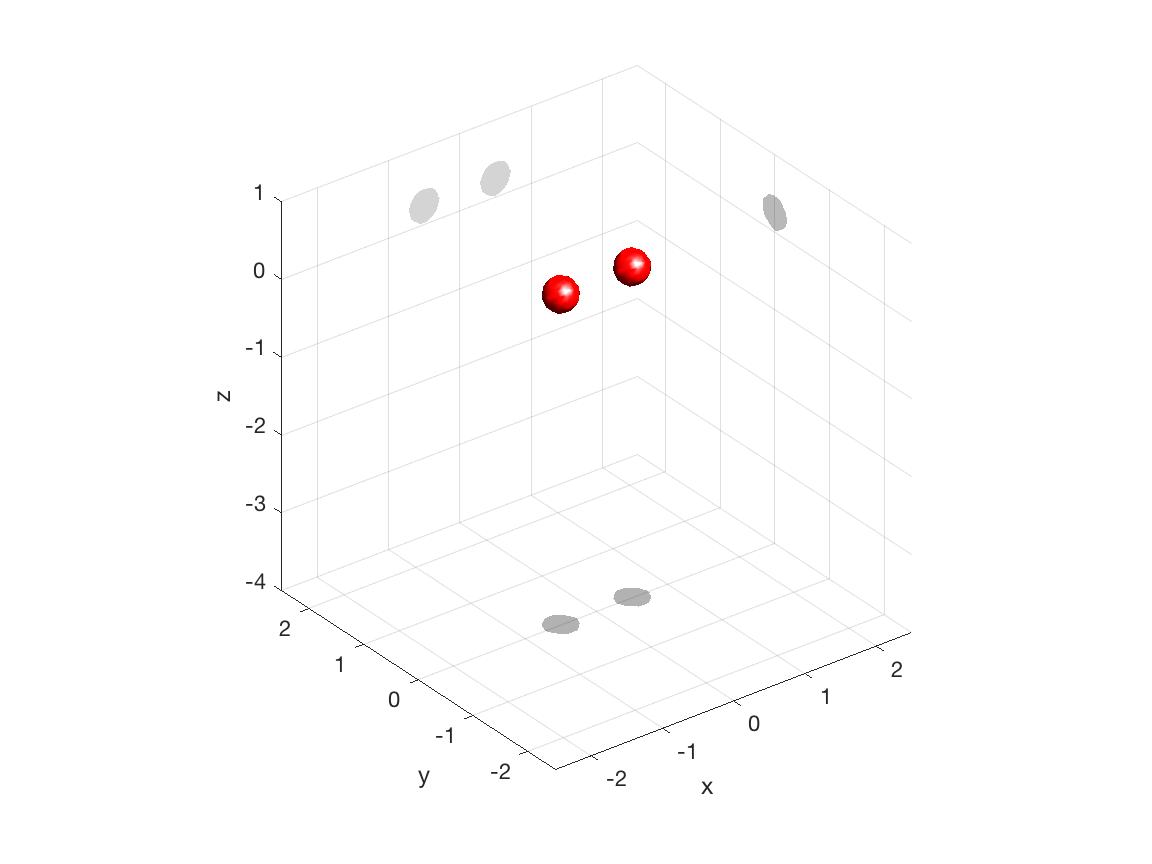}}\\
\subfloat[Reconstruction (cross-sectional
view)]{\includegraphics[width=0.4\textwidth]{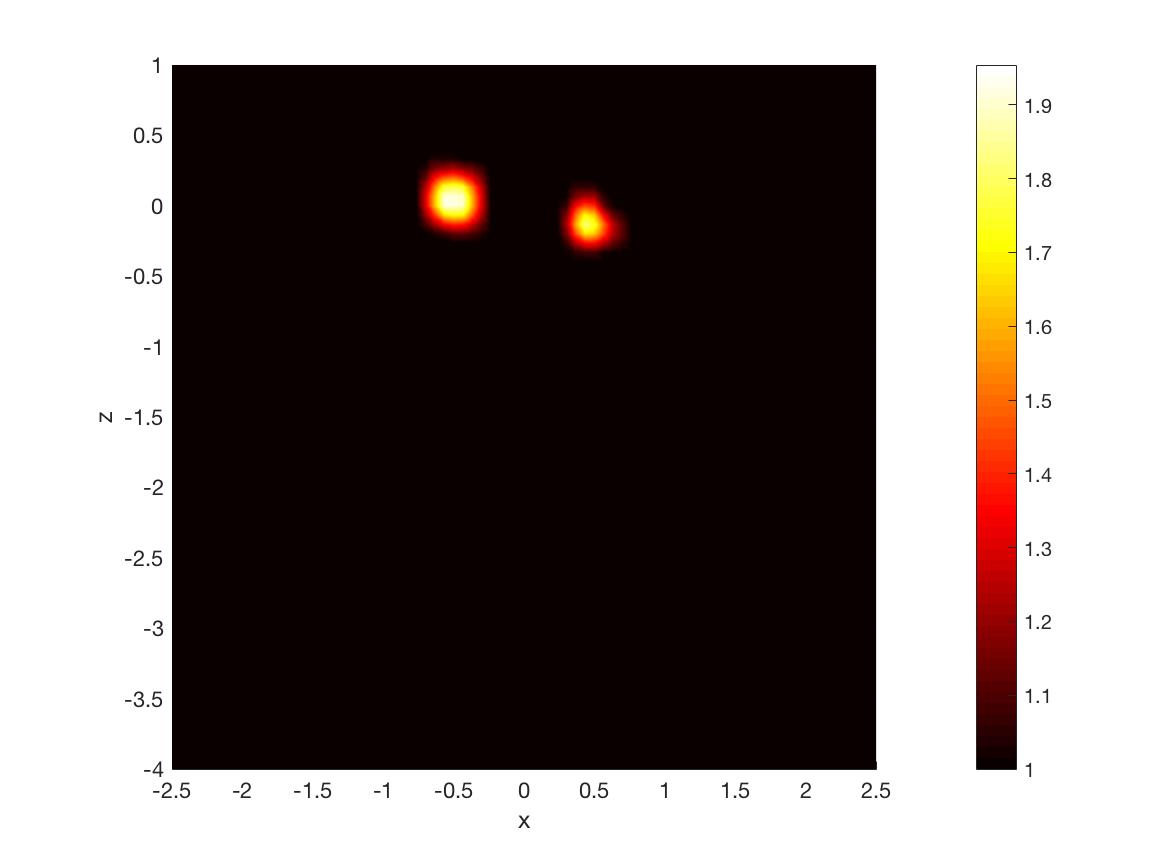}} 
\subfloat[Reconstruction (3D
view)]{\includegraphics[width=0.4\textwidth]{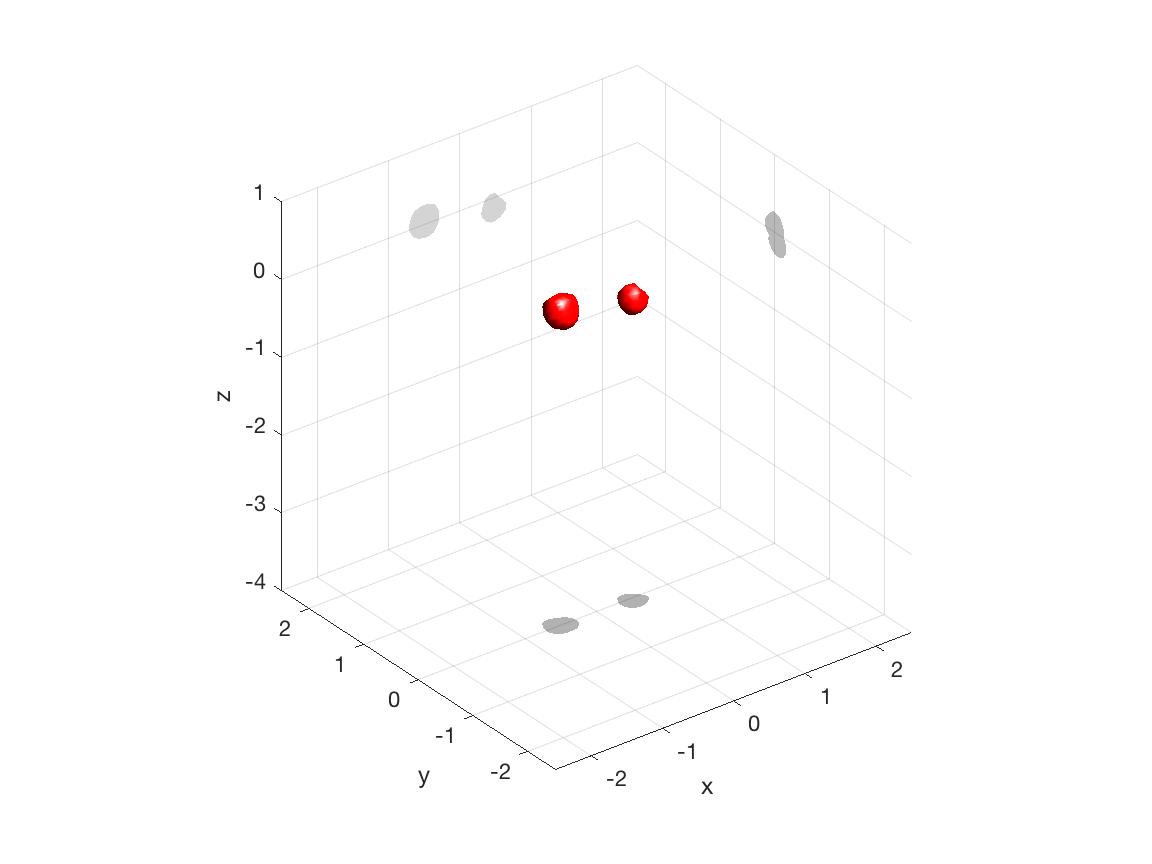}}
\caption{\textit{Reconstruction results for two spherical inclusions: the above item 2.
The maximal value of $c^*(\mathbf{x})$ in both inclusions are 2. Maximal value
of the computed $c_{comp}(\mathbf{x})$ is 1.95 in both inclusions. Hence, the
error in computing this value is 2.5\%}}
\label{fig:res_obj2}
\end{figure}

\begin{figure}[h]
\centering
\subfloat[Exact profile (cross-sectional
view)]{\includegraphics[width=0.4\textwidth]{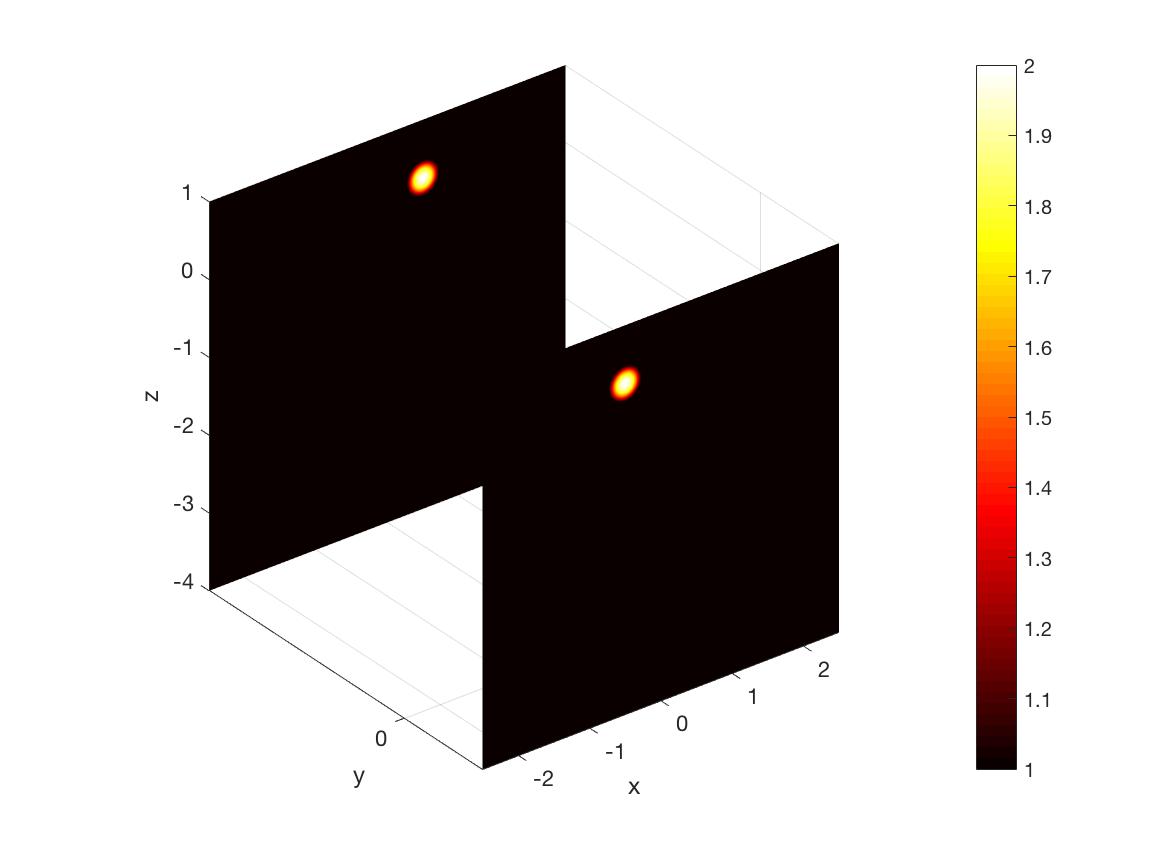}}%
\subfloat[Exact profile (3D
view)]{\includegraphics[width=0.4\textwidth]{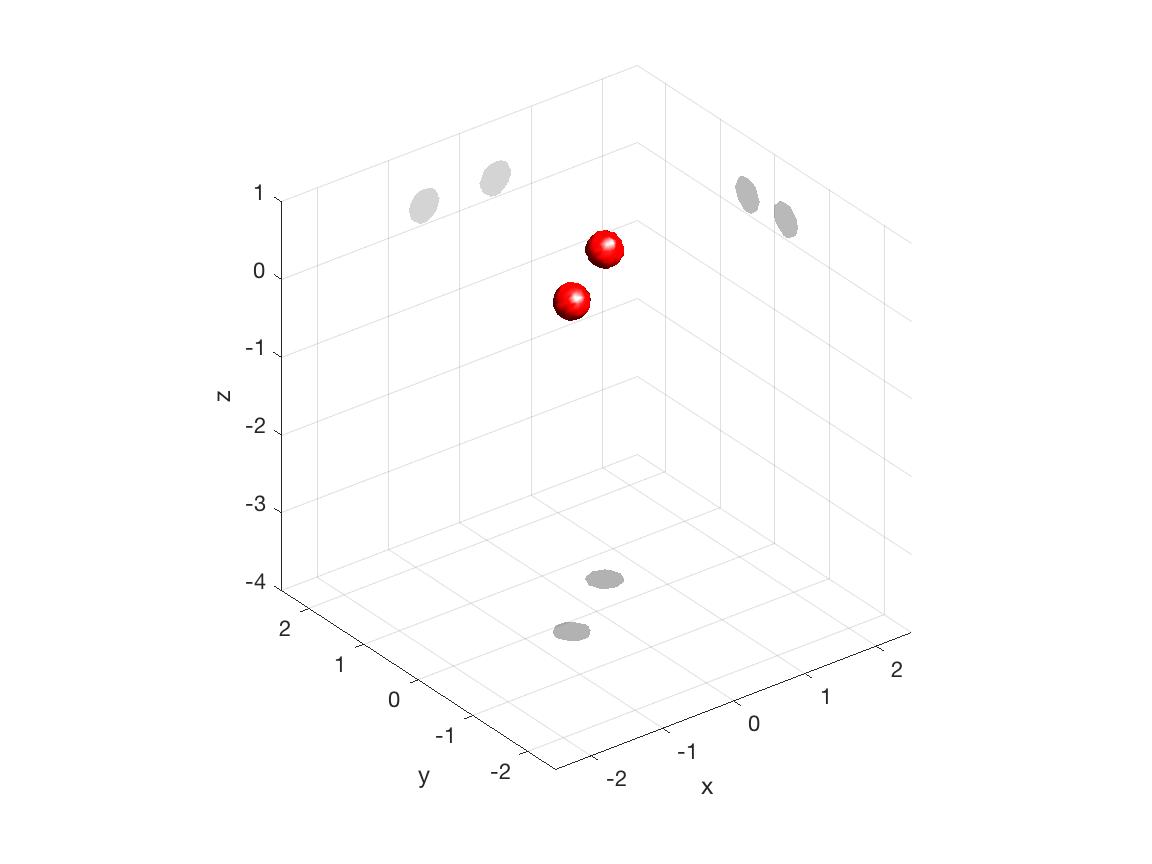}} \\
\subfloat[Reconstruction (cross-sectional
view)]{\includegraphics[width=0.4\textwidth]{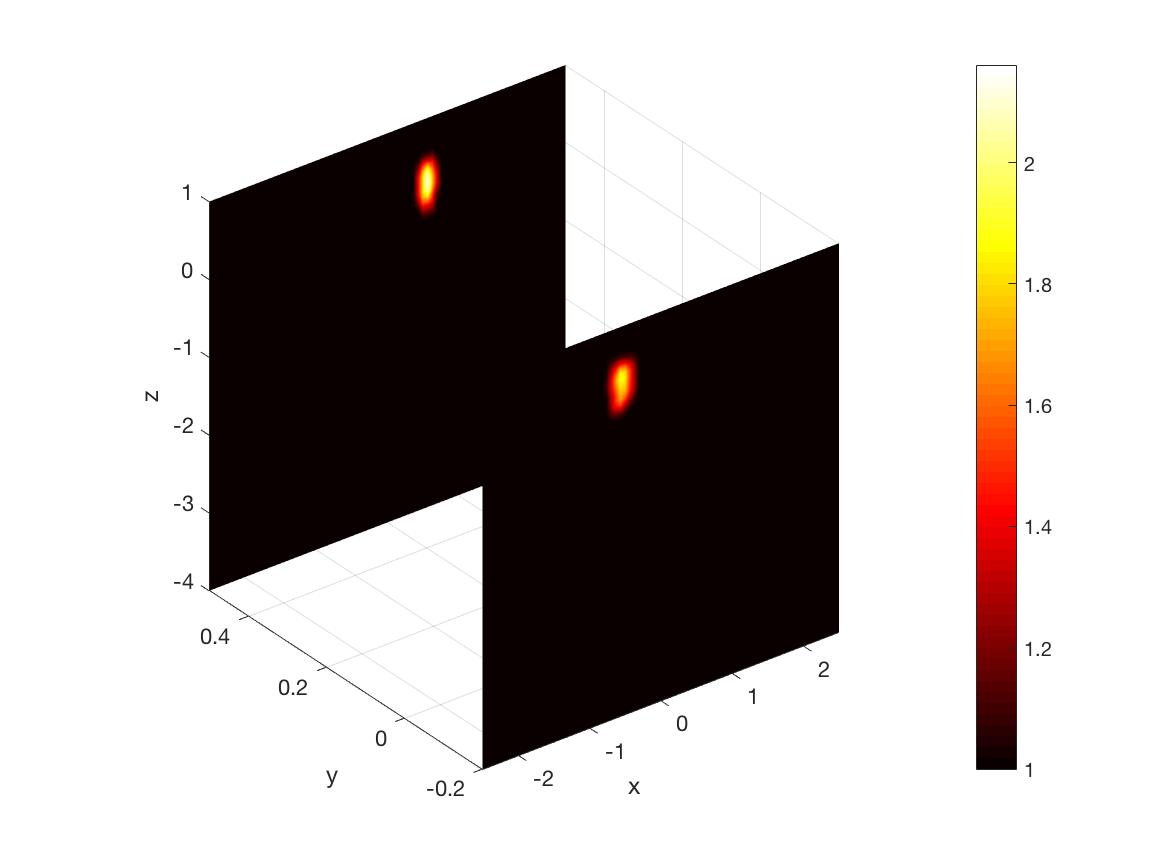}} 
\subfloat[Reconstruction (3D
view)]{\includegraphics[width=0.4\textwidth]{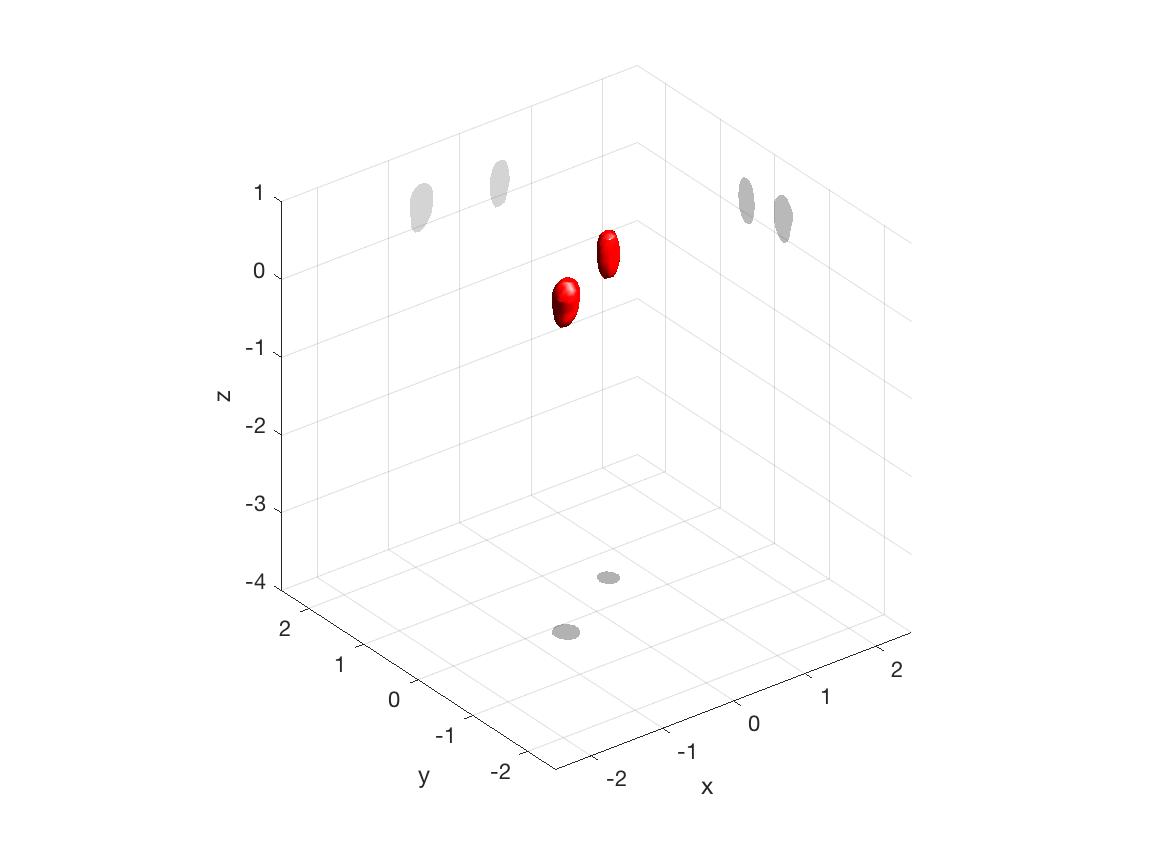}}
\caption{\it Reconstruction results for two non symmetric spherical inclusions: the above item 3. The maximal
value of $c^{\ast }(\mathbf{x})$ in both inclusions is 2. The maximal value of
the computed coefficient $c_{\text{comp}}\left( \mathbf{x}\right) $ is 1.91
and 2.16 in the left and right inclusions respectively. Hence, the error in
computing this value is 4.5\% and 8\% respectively.}
\label{fig:res_obj3}
\end{figure}

%\begin{center}
%\textbf{Liem, Loc,}
%\end{center}
%
%\textbf{In labels for Figures 4-6:} \textcolor{red}{Do you mean 3-5?}
%
%1. Figure 4 \textquotedblleft \emph{one sphere}" should be \textquotedblleft 
%\emph{one spherical inclusion: the above item 1}". Also, \textquotedblleft 
%\emph{this sphere}" should be \textquotedblleft \emph{this inclusion}".
%Otherwise one might think that this is just a surface of a ball. \textcolor{red}{Done}
%
%2. Figure 5 \textquotedblleft \emph{two spheres: the above Case 2}" should
%be \textquotedblleft \emph{two} \emph{spherical inclusions: the above item 2}%
%". Also, \textquotedblleft \emph{spheres}" should be \textquotedblleft \emph{%
%inclusions}". \textcolor{red}{Done}
%
%3. Figure 6.\ Label should be in \emph{italics}. \textquotedblleft non
%symmetric spheres" should be \textquotedblleft non-symmetric \emph{spherical
%inclusions: the above item 3}". Also \textquotedblleft spheres" should be
%\textquotedblleft \emph{inclusions}". \textcolor{red}{Done}

\section{Summary}

\label{sec: 8}

In this paper, we have developed a numerical method for solving a 3D
phaseless inverse scattering problem. Unlike the previous work \cite%
{KlibanovLocKejia:apnum2016} of this group, where overdetermined data were
used and a linearization of the travel time information was applied, we work
here with the data generated by a single measurement event using a single
direction of the incident plane wave and an interval of frequencies. Our
procedure consists of two stages and it does not use any linearization. On
the first stage we reconstruct the first term of the asymptotic expansion at 
$k\rightarrow \infty $ of the function $u_{\text{sc}}\left( \mathbf{x}%
,k\right) $ for $\mathbf{x\in }P_{\text{meas}}.$ As a result, we obtain a
phased coefficient inverse scattering problem.

On the second stage, we solve the latter problem by the globally convergent
numerical method of \cite{Kliba2016}. Our results demonstrate a good
reconstruction accuracy of locations of abnormalities. In addition, the
relative errors in abnormality/background contrasts, which are the maximal
values of the computed coefficients $c_{\text{comp}}\left( \mathbf{x}\right) 
$ of abnormalities, do not exceed 8\% in all cases. Given a significant
complexity of the problem under consideration and 5\% random noise in the
data, we consider this as a quite accurate result.

As to the theoretical part, we prove here uniqueness Theorem 3.2 which
claims that the first term of the asymptotic expansion at $k\rightarrow
\infty $ of the function $u_{\text{sc}}\left( \mathbf{x},k\right) $ can be
uniquely reconstructed from our phaseless data.

\section{Appendix - The full Maxwell's equations and the scalar 3D Helmholtz
equation}

\label{sec: Maxwell}

In this section we present some numerical simulations, which aim to show
that, under some assumptions, the solution of the scalar 3D Helmholtz
equation can be used to approximate such a component of the electric field
satisfying the Maxwell's equations, which is originally incident upon the
medium. We consider two cases here: the backscatter data and the forward
scatter data.

For the backscatter data, we want to verify the use of the Helmholtz
equation in our paper~\cite{Kliba2016} as well as in our works on
experimental microwave data \cite%
{AlekKlibanovLocLiemThanh:anm2017,LiemKlibanovLocAlekFiddyHui:jcp2017,nguyen2017:iip2017}%
. Therefore, we follow the same setup as in~\cite{Kliba2016}. As to the
forward scatter data, we test them to verify the model problem studied in
this paper.

Assume that the scattering objects are isotropic, non-magnetic and that they
are characterized by the dielectric constant $c(\mathbf{x})$, which is a
bounded real-valued function satisfying 
\begin{equation*}
c(\mathbf{x})\geq 1 \mbox{
for all }\mathbf{x}\in \mathbb{R}^{3}\quad \mbox{ and }c(\mathbf{x})=1%
\mbox{
for all }\mathbf{x}\in \mathbb{R}^{3}\setminus \Omega .
\end{equation*}
Let $\mathbf{E}$ be the electric field. The scattering of light in the
frequency domain can be described by the Maxwell's equations for the
electric field as follows: 
\begin{align}  \label{Maxwell}
& \nabla \times \nabla \times \mathbf{E}-k^{2}c(\mathbf{x})\mathbf{E}%
=0,\quad \mathbf{x}\in \mathbb{R}^{3}, \\
& \mathbf{E}(\mathbf{x},k)=\mathbf{E}^{\mathrm{inc}}(\mathbf{x},k)+\mathbf{E}%
^{\mathrm{sc}}(\mathbf{x},k), \\
& \lim_{|\mathbf{x}|\rightarrow \infty }|\mathbf{x}|\left( \nabla \times 
\mathbf{E}^{\mathrm{sc}}\times \hat{\mathbf{x}}-ik\mathbf{E}^{\mathrm{sc}%
}\right) =0,\quad \hat{\mathbf{x}}=\mathbf{x}/|\mathbf{x}|.
\label{radiation}
\end{align}%
The total electric field $\mathbf{E}$ is the sum of the scattered field $%
\mathbf{E}^{\mathrm{sc}}=(E_{1}^{\mathrm{sc}},E_{2}^{\mathrm{sc}},E_{3}^{%
\mathrm{sc}})$ and the incident field $\mathbf{E}^{\mathrm{inc}}$. We note
that the scattered field $\mathbf{E}^{\mathrm{sc}}$ satisfies the
Silver-Muller radiation condition~\eqref{radiation}, which guarantees that
it is an outgoing wave.

We now consider the scattering problem for the scalar Helmholtz equation as
in \cite{Kliba2016}, where the incident wave is a plane wave propagating
along the $z-$direction 
\begin{align}
& \Delta u+k^{2}c(\mathbf{x})u=0,\quad \mathbf{x}\in \mathbb{R}^{3},
\label{Helm} \\
& u=e^{ikz}+u^{\mathrm{sc}}, \\
& \lim_{r\rightarrow \infty }r\left( \partial _{r}u^{\mathrm{sc}}-iku^{%
\mathrm{sc}}\right) =0,\quad r=|\mathbf{x}|.  \label{radiation2}
\end{align}

The measurement square for the backscatter data is $\left\{ \mathbf{x}%
:\left\vert x_{1}\right\vert ,\left\vert x_{2}\right\vert \leq
5, x_{3}=-10\right\} $, while the face of the scatterer is at $\{z=0\}$. For
the forward scatter data, we measure at the same square $P_{\mathrm{meas}}$
as in Section \ref{sect:details}: $P_{\mathrm{meas}}=\{\mathbf{x}%
:|x_{1}|,|x_{2}|\leq 5,x_{3}=2.5\}.$

We observed in our numerical simulation that if the incident field in~%
\eqref{Maxwell}--\eqref{radiation} $\mathbf{E}^{\mathrm{inc}}(\mathbf{x}%
,k)=(0,1,0)e^{ikx_{3}}$, then the second component of the electric scattered
field $E_{2}^{\text{sc}}$ on the measurement square can be well-approximated
by the scattered field of the Helmholtz problem $u_{\text{sc}}$ divided by a
scalar multiplier $d(k)$, defined by 
\begin{equation*}
d(k)=\frac{\max (|u_{\text{sc}}(\mathbf{x},k)|)}{\max (|E_{2}^{\text{sc}}(%
\mathbf{x},k)|)},
\end{equation*}%
where the maximal values are taken on the measurement square. 
%\begin{center}
%\textbf{Liem,}
%\end{center}
%
%\textbf{\textquotedblleft measurement plane" or }$P_{\text{meas}}?$\textbf{\
%Please specify. if plane, then say \textquotedblleft measurement plane }$P$%
%\textbf{" }

The observation above means that the experimental scattering data studied in~%
\cite{Kliba2016, LiemKlibanovLocAlekFiddyHui:jcp2017}, which are supposed to
be the second component of the electric scattered field, can be calibrated
and approximated by the solution to the scattering problem for the scalar
Helmholtz equation. The solution of Maxwell problem~\eqref{Maxwell}--%
\eqref{radiation} and the scalar problem~(\ref{Helm})--(\ref{radiation2})
was computed using the numerical solvers developed in~\cite{Lechl2014} and~%
\cite{Nguye2014a}, respectively.

\begin{figure}[h!]
\centering
\subfloat[]{\includegraphics[width=0.35
\textwidth]{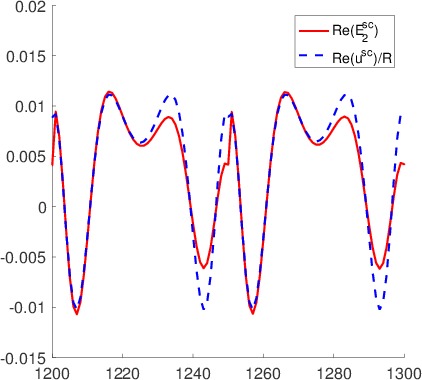}}\hspace{0.5cm} 
\subfloat[]{\includegraphics[width=0.35
\textwidth]{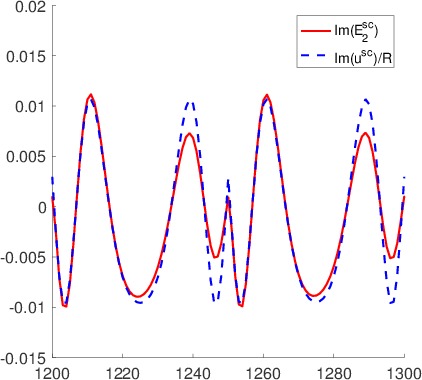}} \\
\subfloat[]{\includegraphics[width=0.35
\textwidth]{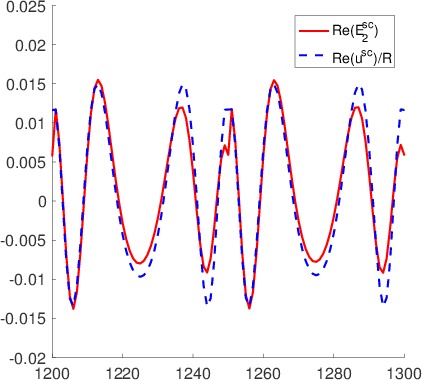}}\hspace{0.5cm} 
\subfloat[]{\includegraphics[width=0.35
\textwidth]{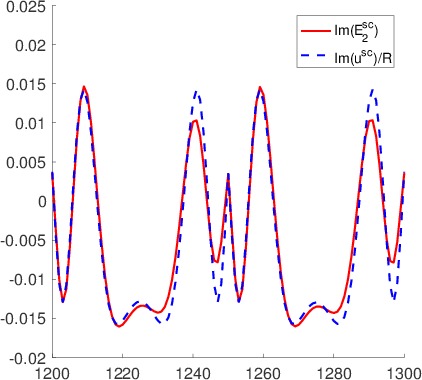}}
\caption{\it The real and imaginary parts of $E^{\mathrm{sc}}_2$ and $u_{\mathrm{%
sc}}/R$ at the measurement points from 1200 to 1300. The scattering object
is a sphere characterized by a smoothly decaying function. (a) and (b) are
for $k = 6.5$, (c) and (d) are for $k = 7.5$.}
\label{fig:exp1}
\end{figure}

\begin{figure}[h!]
\centering
\subfloat[]{\includegraphics[width=0.35
\textwidth]{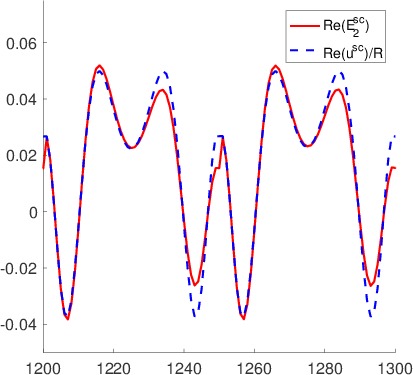}}\hspace{0.5cm} 
\subfloat[]{\includegraphics[width=0.35
\textwidth]{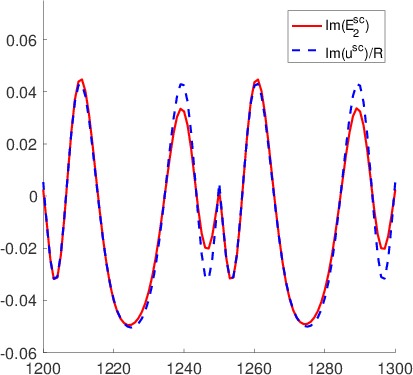}} \\
\subfloat[]{\includegraphics[width=0.35
\textwidth]{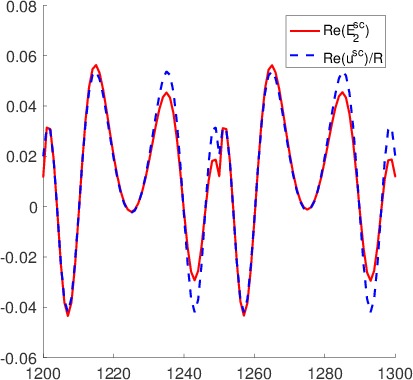}}\hspace{0.5cm} 
\subfloat[]{\includegraphics[width=0.35
\textwidth]{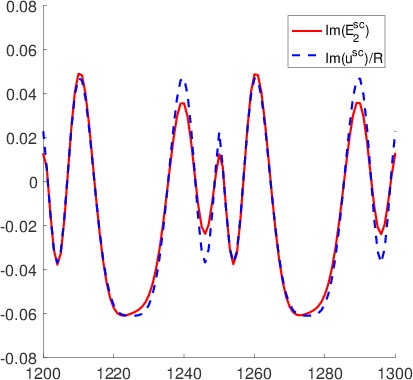}}
\caption{\it The real and imaginary parts of $E^{\mathrm{sc}}_2$ and $u_{\mathrm{%
sc}}/R$ at the measurement points from 1200 to 1300. The scattering object
is a cube. The coefficient $c(\mathbf{x})$ equals 4.5 inside the scatterer
and one elsewhere. (a) and (b) are for $k = 6.5$, (c) and (d) are for $k =
7.5$.}
\label{fig:exp2}
\end{figure}

Now we present four numerical examples. In the first two examples (Figures~%
\ref{fig:exp1} and \ref{fig:exp2}), we consider the backscatter data and the
setup in~\cite{Kliba2016}. More precisely, we consider wave numbers $k=6.5$
and $k=7.5$, the measurement square $[-5,5]^{2}$ is uniformly discretized by 
$50^{2}$ points. To show the fitting of the data in our simulations, we
first transform the $50\times 50$ data matrix in a vector of 2500 points
(using the command : (colon) in MATLAB). Then we choose the data points from
1200 to 1300, where the signals of the scattered fields are the strongest,
to present them in Figures~\ref{fig:exp1} and \ref{fig:exp2}. We note that
the results are similar for other measurement points, where the scattered
fields are weaker. Figure~\ref{fig:exp1} is dedicated to the case of a
spherical scattering object represented by a smooth function $c(\mathbf{x})$%
, which is similar to the model in Section~\ref{sect:numerical} except $\max
\{c(\mathbf{x})\}=4.5$ instead of 2. We consider in Figure~\ref{fig:exp2}
the case of a rectangular scattering object represented by a piecewise
constant function $c(\mathbf{x})$ with jumps across the boundary of the
scatterer ($c(\mathbf{x})=4.5$ inside the scatterer).

\begin{figure}[h!]
\centering
\subfloat[One smooth inclusions, k =
20.5]{\includegraphics[width=0.35\textwidth]{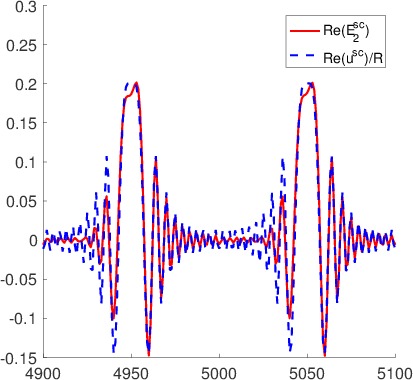}}%
\hspace{0.5cm} 
\subfloat[One smooth
inclusions, k =
20.5]{\includegraphics[width=0.35\textwidth]{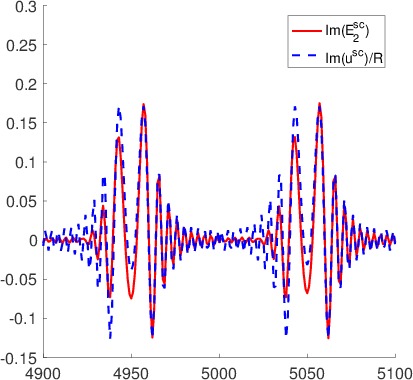}}%
\\
\subfloat[One smooth inclusions, k =
21.5]{\includegraphics[width=0.35\textwidth]{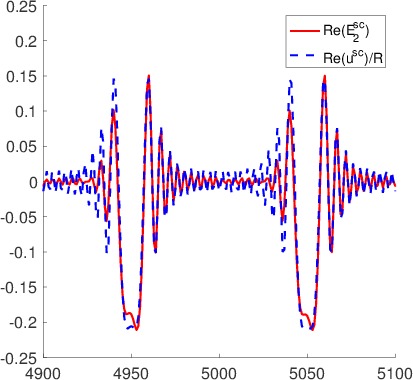}}%
\hspace{0.5cm} 
\subfloat[One smooth
inclusions, k =
21.5]{\includegraphics[width=0.35\textwidth]{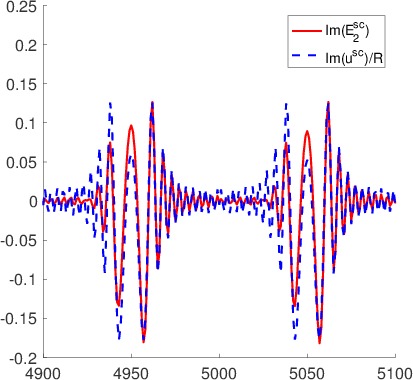}}
\caption{\it The real and imaginary parts of $E^{\mathrm{sc}}_2$ and $u_{\mathrm{%
sc}}/R$ at the measurement points from 4900 to 5100. The model for the
scattering object, which is a sphere, is given in the case of one inclusion
in Section~\protect\ref{sect:numerical}. (a) and (b) are for $k = 20.5$, (c)
and (d) are for $k = 21.5$.}
\label{fig:exp3}
\end{figure}

\begin{figure}[h]
\centering
\subfloat[Two smooth inclusions, k =
20.5]{\includegraphics[width=0.35\textwidth]{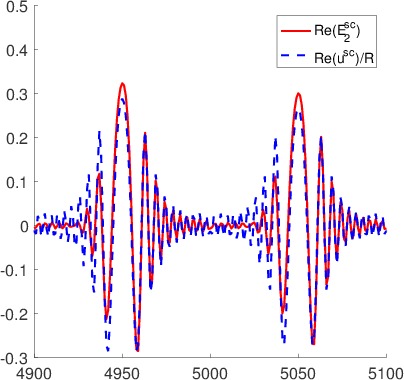}}
\hspace{0.5cm} 
\subfloat[Two smooth inclusions, k =
20.5]{\includegraphics[width=0.35\textwidth]{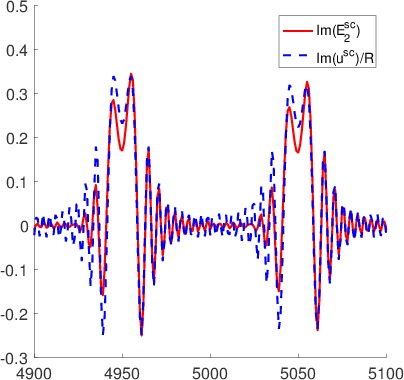}}
\\
\subfloat[Two smooth inclusions, k =
21.5]{\includegraphics[width=0.35\textwidth]{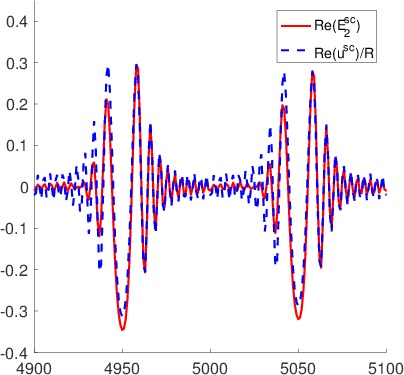}}%
\hspace{0.5cm} 
\subfloat[Two smooth inclusions, k =
21.5]{\includegraphics[width=0.35\textwidth]{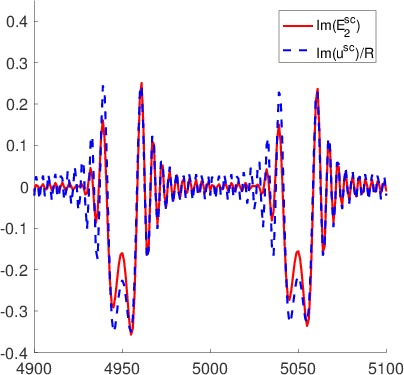}}
\caption{\it The real and imaginary parts of $E_{2}^{\mathrm{sc}}$ and $u_{%
\mathrm{sc}}/R$ at the measurement points from 4900 to 5100. The model for
the scattering object, which is two spheres, is given in the case of two
inclusions in Section~\protect\ref{sect:numerical}. (a) and (b) are for $%
k=20.5$, (c) and (d) are for $k=21.5$.}
\label{fig:exp4}
\end{figure}

In Figure~\ref{fig:exp3} and Figure~\ref{fig:exp4} we present the numerical
simulations for the setup considered in Section~\ref{sect:numerical}. More
precisely, we consider the case of one inclusion and the case of two
inclusions there, where wave number $k$ is 20.5 and 21.5. We have $100^{2}$
measurement points and choose to present the points from 4900 to 5100, where
scattered fields seem to be strongest.

\begin{center}
\textbf{Acknowledgements}
\end{center}

This work was supported by the Office of Naval Research grant
N00014-15-1-2330 as well as by the US Army Research Laboratory and US Army
Research Office grant W911NF-15-1-0233. In addition, the effort of L.H.
Nguyen was partially supported by research funds FRG 111172 provided by
University of North Carolina at Charlotte. The authors are grateful to
Professor Vasilii Astratov from Center for Optoelectronics and Optical
Communications of University of North Carolina at Charlotte for many
fruitful discussions.

\providecommand{\noopsort}[1]{}

\end{document}